\documentclass[11pt]{article}
\usepackage{bbm}
\usepackage{mathrsfs}
\usepackage{amsfonts}
\usepackage{amssymb}
\usepackage{amssymb,amsmath,amsthm, amsfonts}
\usepackage{indentfirst}
\usepackage{titlesec}
\usepackage{geometry}
\usepackage{cite}

\makeatletter
\@addtoreset{equation}{section}
\newtheorem{theorem}{Theorem}[section]
\newtheorem{corollary}[theorem]{Corollary}
\newtheorem{lemma}[theorem]{Lemma}
\newtheorem{remark}[theorem]{Remark}
\newtheorem{definition}[theorem]{Definition}
\newtheorem{prop}[theorem]{Proposition}

\makeatother

\titleformat*{\section}{\large\bfseries}
\titleformat*{\subsection}{\normalsize\bfseries}
\begin{document}

\title{Mean field stochastic control under sublinear expectation\footnotemark[1]}
 \author{Rainer Buckdahn$^{1,3}$,\,\, Bowen He$^{2,**}$,\,\, Juan Li$^{2,3,**}$ \\
 {$^1$\small Laboratoire de Math\'{e}matiques de Bretagne Atlantique, Univ Brest,}\\
	{\small UMR CNRS 6205, 6 avenue Le Gorgeu, 29200 Brest, France.}\\
{$^2$\small  School of Mathematics and Statistics, Shandong University, Weihai,}
	{\small Weihai 264209, P. R. China.}\\
	{$^3$\small  Research Center for Mathematics and Interdisciplinary Sciences, Shandong University,}\\
	{\small Qingdao 266237, P. R. China.}\\
 {\small{\it E-mails: rainer.buckdahn@univ-brest.fr,\,\ hbw@mail.sdu.edu.cn,\,\ juanli@sdu.edu.cn.}}
 \date{November 09, 2022}}
\renewcommand{\thefootnote}{\fnsymbol{footnote}}
\footnotetext[1]{The work is supported by the NSF of P.R. China (NOs. 12031009, 11871037), National Key R and D Program of China (NO. 2018YFA0703900), and NSFC-RS (No. 11661130148; NA150344).\\
\ \ ${    }^{**}$ Corresponding authors.}
\maketitle
\noindent\textbf{Abstract.}
Our work is devoted to the study of Pontryagin's stochastic maximum principle for a mean-field optimal control problem under Peng's $G$-expectation. The dynamics of the controlled state process is given by a stochastic differential equation driven by a $G$-Brownian motion, whose coefficients depend not only on the control, the controlled state process but also on its law under the $G$-expectation.
Also the associated cost functional is of mean-field type. Under the assumption of a convex control state space we study the stochastic maximum principle, which gives a necessary optimality condition for control processes. Under additional convexity assumptions on the Hamiltonian it is shown that this necessary condition is also a sufficient one. The main difficulty which we have to overcome in our work consists in the differentiation of the $G$-expectation of parameterized random variables. As particularly delicate it turns out to handle with the $G$-expectation of a function of the controlled state process inside the running cost of the cost function. For this we have to study a measurable selection theorem for set-valued functions whose values are subsets of the representing set of probability measures for the $G$-expectation.

\vskip 0.5cm



\vskip 1cm
\newpage
\section{Introduction}
Our work brings together two important subjects of actual intensive research, mean-field problems popularised by Lasry and Lions' pioneering work \cite{18} on mean-field games in 2007  on one side and stochastic control under Peng's sublinear $G$-expectation (see, e.g., \cite{9} and \cite{10}) on the other side. More precisely, we study Pontryagin's stochastic maximum principle (SMP) for a stochastic control problem over a $G$-expectation space, whose dynamics are given by a controlled $G$-stochastic differential equation ($G$-SDE) whose coefficients do not only depend on the control process and the associated controlled state process but also on its law under the $G$-expectation, which we consider as the $G$-expectation of a function of the controlled state space. Also in the associated cost functional both the terminal cost function and the running cost function depend not only on the controlled state process but also on its law with respect to (w.r.t.) the $G$-expectation.

Mean-field SDEs in form of McKean-Vlasov equations have
been studied for a long time and have found a lot of applications in different domains. Recently, with their seminal paper \cite{18} on mean-field games and their applications in economics, finance and game theory, Lasry and Lions have given new impulses to this research topic, opened the way to new applications and attracted a lot of researchers to this topic.
One of these applications is the study of mean-field stochastic optimal control problems. Motivated by the rich literature on the stochastic maximum principle in the classical stochastic control, for example, Peng's SMP \cite{4}, different authors studies the stochastic maximum principle in the context of mean-field control problems. Let us namely mention the work by Buckdahn, Djehiche and Li \cite{5} in 2011, where the coefficients of the mean-field SDEs depend on the solution process, its expectation and the control. Li \cite{7} studied Pontryagin's SMP for mean-field SDEs, and obtained  necessary and sufficient conditions for the optimality of a control process, while Buckdahn, Li, Ma \cite{8,17} studied the optimal control problem for a class of general mean-field SDEs, in which the coefficients depend non linearly on both the state process as well as on its law. They extended the SMP of Buckdahn et al. \cite{5} to this general case.  Acciaio et al. \cite{20} studied the stochastic maximum principle for an extended mean-field control problem.

However, for instance, in economy and in finance a vast field of applications requires to model Knightian uncertainty. Inspired by financial problems with uncertainty, Peng \cite{9} introduced a fully non linear expectation, called $G$-expectation $\hat{\mathbb{E}}[\cdot]$, and he proved that it can well characterize the Knightian uncertainty. Under this $G$-expectation framework a new type of Brownian motion, the so-called $G$-Brownian motion, has been introduced and the stochastic calculus with respect to the $G$-Brownian motion has been developed. Recently, Hu et al. \cite{12,13} developed the SDE and BSDE theory in this $G$-expectation framework. And they also studied the SMP for stochastic optimal control problems under $G$-expectation or uncertainty (see \cite{14,15}).

For the reasons explained above we study the SMP for a mean-field stochastic control problem under $G$-expectation. We consider a stochastic control problem where the dynamic of the state process is given by a
stochastic differential equation driven by a $G$-Brownian motion ($G$-SDE) of mean-field type. That is, the coefficients do not only depend on the control and the controlled state process but also on some functional
of the law of the state process under the sublinear expectation.  More precisely, we consider the dynamics
$$\left\{\begin{array}{l}
d X_{t}^{u}=\sigma(X_{t}^{u}, \hat{\mathbb{E}}[\varphi_{1}(X_{t}^{u})], u_{t}) d B_{t}+b(X_{t}^{u} , \hat{\mathbb{E}}[\varphi_{2}(X_{t}^{u})], u_{t}) d t
+\beta(X_{t}^{u}, \hat{\mathbb{E}}[\varphi_{3}(X_{t}^{u})], u_{t}) d \langle B \rangle_{t},\\
X_{0}^{u}=x_{0}\in \mathbb{R}^n, \quad t \in[0, T],
\end{array} \right.$$
for some functions $b$, $\sigma$, $\beta$, $\varphi_i, ~i=1,2,3$, and  the $G$-Brownian motion $B=(B_t)$. The admissible control process $u=(u_t)$ takes its values in a convex state space $U$.
The objective is to minimize the associated cost functional of the form
$$
J(u):=\hat{\mathbb{E}}\big[\Phi(X_{T}^{u},\hat{\mathbb{E}}[\varphi_{4}(X_{T}^{u})])+\int_{0}^{T} l(t, X_{t}^{u}, \hat{\mathbb{E}}[\varphi_{5}(X_{t}^{u})], u_{t}) d t\big],
$$
for given functions $\Phi$, $l$ and $\varphi_i,~ i=4,5$. Also this cost functional is of mean-field type, as
the functions $\Phi$ and $l$, depend on the law under the sublinear expectation of the state process.

In this paper we derive necessary and sufficient conditions for optimality of this
control problem in form of a stochastic maximum principle for a convex action space, using
the convex perturbation. The stochastic
maximum principle involves solving a family of adjoint equations, backward SDEs (BSDEs).

As concerns previous works related with the SMP under sublinear expectation, we have to mention mainly the recent works by Biagini, Meyer-Brandis and {{\O}}ksendal \cite{27}, Sun \cite{28} and Hu and Ji \cite{14}. In \cite{27} the authors study a stochastic control problem (without mean-field term), composed of a forward $G$-SDE and a cost functional under $G$-expectations (also without mean-field term). The sufficient but also the necessary optimality conditions for a control process $\hat{u}$ they give need the assumption that in the adjoint equation, a $G$-BSDE (see Definition \eqref{def2.10}),  the non increasing $G$-martingale $K$ is identically equal to zero. In \cite{28} Sun studies a controlled system of $G$-forward and $G$-backward SDEs with solution $(X^{u},Y^{u},Z^{u})$, and he associates the cost functional $J(u)=\hat{\mathbb{E}}[\psi(u)]$ with $\psi(u)=\phi(X_T^{u})+\int_0^Tl(t,X_t^{u},Y_t^{u},Z_t^{u},u_t)dt+\gamma(Y_0^{u})$. In the deduction of the sufficient optimality condition for a control $\hat{u}$ he uses convexity assumptions. Also in a non mean-field context, Hu and Ji \cite{14}  study a system of forward and backward $G$-SDEs, they consider as cost functional $J(u)=Y^{u}_0,$ and they investigate the SMP. For this they show namely that the cost functional for the perturbed optimal control $\lambda\mapsto J(\hat{u}+\lambda(u-\hat{u}))$ is right-differentiable at $\lambda=0$, and they use the special form of this derivative and an application of Sion's minimax theorem to derive a necessary optimality condition of the optimal control $\hat{u}$. Their approach depends on the related $G$-BSDEs.

Inspired by above works we study the mean-field stochastic control problem under $G$-expectation. We investigate the SMP and give a necessary optimality condition for the optimal control and also a sufficient one for the optimality of a control. However, the fact that we have to do not only with the $G$-expectation of the definition of the cost functional $J(u)$ but also with the $G$-expectations $\hat{\mathbb{E}}[\varphi_i(X_t^{u})],\, 1\le i\le 5$, involves new difficulties. So, for instance, in the general case, Sion's minimax theorem cannot be applied.
As it plays a crucial role, a whole section (Section 4) is devoted to the study of the derivative of functions of laws under $G$-expectation. For this we begin with the easy observation that, given two random variables $\xi,\, \eta$, the function $\lambda\mapsto F(\lambda)=\hat{\mathbb{E}}[\xi+\lambda\eta]$ is convex, i.e., the right but also the left derivatives $F'_{+}(\lambda)$ and $F'_{-}(\lambda)$, respectively, exist. We determine them in a more direct approach than that in \cite{14}, without passing through the associated $G$-BSDE, and we also associate some essential result which has its own interest (see Proposition \ref{prop4.2}). The results are extended to the derivative of functions $F_f(\xi):=\displaystyle\sup_{P\in\mathcal{P}}f(P_\xi)$, where $\mathcal{P}$ represents $\hat{\mathbb{E}}[\cdot]$ (see Theorem \ref{th2.7}). As the derivative of this latter function is not directly used for our SMP approach but has its own interest, it is shifted to Appendix 1. Section 5 is devoted to deduce the necessary optimality condition for an optimal control. Our main result is Theorem \ref{th5.5}.
The main difficulty here stems from the fact that our coefficients depend also on $\hat{\mathbb{E}}[\varphi_i(X_t^{u})],\, 1\le i\le 5$, and so all their derivatives have to be considered. The most delicate part comes from the dependence of the running cost $l$ on $\hat{\mathbb{E}}[\varphi_5(X_t^{u})]$.  To handle these difficulties we need a measurable selection theorem for a mapping $[0,T]\ni\mapsto\mathcal{P}_{\xi_t|\eta_t}\subset \mathcal{P}$ (see Theorem \ref{thA.2}). Proving that this mapping is a weakly measurable set-valued function, we can use the Kuratowski and Ryll-Nardzewski measurable selection theorem to get Theorem \ref{thA.2}; see Appendix 2.  For the case that the running cost coefficient $l$ does not depend $\hat{\mathbb{E}}[\varphi_5(X_t^{u})]$, Sion's minimax theorem can be used to simplify the necessary optimality condition considerably; see Theorem \ref{th5.6}. The second part of Section 5 is devoted to the study of a sufficient optimality condition for the general case and to an example.

Our paper is organized as follows: In Section 2, we recall some basic notions of
$G$-expectation and results of $G$-SDEs and $G$-BSDEs. Section 3 introduces the formulation of the mean-field stochastic control problem, and Section 4 is devoted to the study of the derivative of the $G$-expectation of parameter depending random variables. In Section 5 we study the SMP and the Appendix is devoted to an extension of the studies made in Section 4 and to the proof of our measurable selection theorem.

\section{Preliminaries}
In this section, we review some notations and results in the $G$-expectation framework, which
are mainly concerned with the $G$-It$\rm\hat{o}$ calculus and BSDEs driven by a $G$-Brownian motion.  More relevant
details can be found in \cite{9,10,11,12,19}.
\subsection{$G$-expectation space}

 Let $\Omega$ be a given non empty set and $\mathcal{H}$ be a linear space of real-valued functions on $\Omega$ such that, for all $d\ge 1$, if $X_{1}, \ldots, X_{d} \in \mathcal{H}$, then also $\varphi(X_{1}, X_{2}, \ldots, X_{d}) \in \mathcal{H}$ for every $\varphi \in C_{b . L i p}(\mathbb{R}^{d})$, where
$C_{b . L i p}(\mathbb{R}^{d})$ is the space of bounded Lipschitz functions on $\mathbb{R}^{d}$. The set $\mathcal{H}$ is considered as the space of random variables.
\begin{definition}\rm
A sublinear expectation $\hat{\mathbb{E}}$ on $\mathcal{H}$ is a functional $\hat{\mathbb{E}}: \mathcal{H} \rightarrow \mathbb{R}$ having the following properties: For each $X, Y \in \mathcal{H}$,

(i) Monotonicity: $\hat{\mathbb{E}}[X] \geq \hat{\mathbb{E}}[Y]$, if $X \geq Y$;

(ii) Constant preserving: $\hat{\mathbb{E}}[c]=c$, for $c \in \mathbb{R} ;$

(iii) Sub-additivity: $ \hat{\mathbb{E}}[X+Y] \leq \hat{\mathbb{E}}[X]+\hat{\mathbb{E}}[Y];$

(iv) Positive homogeneity: $\hat{\mathbb{E}}[\lambda X]=\lambda \hat{\mathbb{E}}[X]$, for all real $\lambda \geq 0$.

\noindent The triple $(\Omega, \mathcal{H}, \hat{\mathbb{E}})$ is called a sublinear expectation space.
\end{definition}
\begin{definition}\rm
Two $d$-dimensional random vectors $X_{1}$ and $X_{2}$ defined, respectively, on sublinear expectation spaces $(\Omega_{1}, \mathcal{H}_{1}, \hat{\mathbb{E}}_{1})$ and $(\Omega_{2}, \mathcal{H}_{2}, \hat{\mathbb{E}}_{2})$ are called identically distributed, denoted by $X_{1} \stackrel{d}{=} X_{2}$, if
$$
\hat{\mathbb{E}}_{1}[\varphi(X_{1})]=\hat{\mathbb{E}}_{2}[\varphi(X_{2})], \quad \text { for every } \varphi \in C_{b . L i p}(\mathbb{R}^{d}).
$$
\end{definition}
\begin{definition}\rm
On the sublinear expectation space $(\Omega, \mathcal{H}, \hat{\mathbb{E}})$, an $n$-dimensional random vector $Y$ is said to be independent of a $d$-dimensional random vector $X$, denoted by $Y \perp X$, if
$$\hat{\mathbb{E}}[\varphi(X, Y)]=\hat{\mathbb{E}}[\hat{\mathbb{E}}[\varphi(x, Y)]_{x=X}], \quad \text { for every } \varphi \in C_{b . L i p}(\mathbb{R}^{d+n}).$$

A $d$-dimensional random vector $\bar{X}$ is said to be an independent copy of $X$ if $\bar{X} \stackrel{d}{=} X$ and $\bar{X} \perp X$.
\end{definition}
\begin{prop}
Let $X,Y\in\mathcal{H}$ be such that $\hat{\mathbb{E}}[Y]=-\hat{\mathbb{E}}[-Y]$. Then we have
$$
\hat{\mathbb{E}}[X+Y]=\hat{\mathbb{E}}[X]+\hat{\mathbb{E}}[Y].
$$
\end{prop}\rm
\begin{definition}\rm
 A $d$-dimensional random vector $X$ defined on $(\Omega, \mathcal{H}, \hat{\mathbb{E}})$ is called $G$-normally distributed if for any $a, b \geq 0$,
$$
a X+b \bar{X} \stackrel{d}{=} \sqrt{a^{2}+b^{2}} X,
$$
where $\bar{X}$ is an independent copy of $X$. Here the letter $G$ denotes the function $G(A):=$ $\frac{1}{2} \hat{\mathbb{E}}[\langle A X, X\rangle]$, for $A \in \mathbb{S}(d)$, where $\mathbb{S}(d)$ is the space of all $d \times d$ symmetric matrices.
\end{definition}
Throughout this paper, we denote by $\Omega:=C([0, \infty) ; \mathbb{R}^{d})$ the space of all $\mathbb{R}^{d}$-valued continuous paths $(\omega_{t})_{t \geq 0}$, equipped with the distance
$$
\rho_{d}(\omega^{1}, \omega^{2}):=\sum_{i=1}^{\infty} \frac{1}{2^{i}}(\|\omega^{1}-\omega^{2}\|_{C([0, i] ; \mathbb{R}^{d})} \wedge 1),
$$
where $\|\omega^{1}-\omega^{2}\|_{C([0, T] ; \mathbb{R}^{d})}:=\displaystyle\max _{t \in[0, T]}|\omega_{t}^{1}-\omega_{t}^{2}|,$ for $T>0 .$ Given any $T>0$, we also define
$\Omega_{T}:=\{(\omega_{t \wedge T})_{t \geq 0}: \omega \in \Omega\}.$

Let $B_{t}(\omega):=\omega_{t}$, $\omega \in \Omega, t \geq 0$, be the coordinate process on $\Omega$. We introduce the space
\begin{equation*}
\begin{aligned}
L_{i p}(\Omega_{T})\! :=\!\{\varphi(B_{t_{1}}, B_{t_{2}}\!\!-\!B_{t_{1}}, \cdots\!, B_{t_{n}}\!\!-\!B_{t_{n-1}}\!)\!:\! n \in \mathbb{N},
~0 \leq t_{1}<t_{2} \cdots\!<t_{n} \leq T, ~\varphi\! \in\! C_{b . L i p}(\mathbb{R}^{d \times n})\},
\end{aligned}
\end{equation*}
as well as
$\displaystyle
L_{i p}(\Omega):=\bigcup_{m=1}^{\infty} L_{i p}(\Omega_{m}) .
$

The $G$-expectation on $L_{i p}(\Omega)$ is defined by
$$
\hat{\mathbb{E}}[X]:=\widetilde{\mathbb{E}}[\varphi(\sqrt{t_{1}} \xi_{1}, \sqrt{t_{2}-t_{1}} \xi_{2}, \ldots, \sqrt{t_{n}-t_{n-1}} \xi_{n})],
$$
for all $X=\varphi(B_{t_{1}}, B_{t_{2}}-B_{t_{1}}, \ldots, B_{t_{n}}-B_{t_{n-1}}),~ n\ge 1,\, 0 \leq t_{1}<\cdots<t_{n}<\infty$, where $\{\xi_{i}\}_{i=1}^{n}$
is a collection of $n$ $d$-dimensional identically distributed random variables on a sublinear expectation space $(\widetilde{\Omega}, \widetilde{\mathcal{H}}, \widetilde{\mathbb{E}})$ such that, for all $1\le i\le n$, $\xi_{i}$ is $G$-normally distributed and independent of $(\xi_{1}, \ldots, \xi_{i-1})$. Then under $\hat{\mathbb{E}}$, the coordinate process $B_{t}=(B_{t}^{1}, \ldots, B_{t}^{d})$ is a $d$-dimensional $G$-Brownian motion defined by the following properties:

(a) $B_{0}=0$;

(b) For every $t ,\thinspace s \geq 0$, the increment $B_{t+s}-B_{t}$ is independent of $(B_{t_{1}}, \ldots, B_{t_{n}})$, for all $n \in \mathbb{N}$ and $0 \leq t_{1} \leq \cdots \leq t_{n} \leq t$;

(c) $B_{t+s}-B_{t} \stackrel{d}{=} \sqrt{s} \xi$, for $t, s \geq 0$, where $\xi$ is $G$-normally distributed.
\begin{remark}\rm
(i) It is easy to check that the $G$-Brownian motion is symmetric, i.e., $(-B_{t})_{t \geq 0}$ is also a $G$-Brownian motion.

(ii) If, in particular, $G(A)=\frac{1}{2} \operatorname{tr}(A)$, then the $G$-expectation is just a linear expectation with respect to the Wiener measure $P$, i.e., $\hat{\mathbb{E}}=E_{P}$, and the $G$-Brownian motion is a classical Brownian motion over $(\Omega,\mathcal{B}(\Omega),P)$ ($\mathcal{B}(\Omega)$ denotes the Borel $\sigma$-field over $(\Omega,\rho_d)$).
\end{remark}
The conditional $G$-expectation (knowing $\mathcal{B}(\Omega_t)$) for $X=\varphi(B_{t_{1}}, B_{t_{2}}-B_{t_{1}}, \ldots, B_{t_{n}}-B_{t_{n-1}})$ at $t=t_{j}, 1 \leq j \leq n$,
is defined by
$$
\hat{\mathbb{E}}_{t_{j}}[X]:=\phi(B_{t_{1}}, B_{t_{2}}-B_{t_{1}}, \ldots, B_{t_{j}}-B_{t_{j-1}}),
$$
where $\phi(x_{1}, \ldots, x_{j})=\hat{\mathbb{E}}[\varphi(x_{1}, \ldots, x_{j}, B_{t_{j+1}}-B_{t_{j}}, \ldots, B_{t_{n}}-B_{t_{n-1}})] .$

For every $p \geq 1$, we denote by $L_{G}^{p}(\Omega_{t})$ $(L_{G}^{p}(\Omega),$ resp.) the completion of $L_{i p}(\Omega_{t})$ $(L_{i p}(\Omega),$ resp.) under the norm $\|X\|_{p}:=(\hat{\mathbb{E}}[|X|^{p}])^{1 / p} .$ The conditional $G$-expectation $\hat{E}_{t}[\cdot]$ ($t\ge 0$) can be extended continuously to $L_{G}^{1}(\Omega)$.

 We recall the following representation theorem.
\begin{theorem}[\!\!\cite{19,21}]\label{th2.7}
Let

\centerline{$\mathcal{P}=\{P$   probability on $(\Omega, \mathcal{B}(\Omega)): E_{P}[X] \leq \hat{\mathbb{E}}[X]$,
for all $X \in L_{G}^{1}(\Omega)\}$.}

\noindent Then $\mathcal{P}\not=\emptyset$ is a convex, weakly compact subset of the space $\mathcal{P}(\mathbb{R}^d)$ of all probability measures over $(\mathbb{R}^d,\mathcal{B}(\mathbb{R}^d))$ endowed with the topology of weak convergence, and
$$
\hat{\mathbb{E}}[\xi]=\sup _{P \in \mathcal{P}} E_{P}[\xi],  \text { for all } \xi \in L_{G}^{1}(\Omega).
$$
The set $\mathcal{P}$ is said to represent $\hat{\mathbb{E}}$.
\end{theorem}
 The following definition introduces the notion of distributions of random variables under $G$-expectation.
\begin{definition}\rm
Let $X=(X_{1}, \cdots, X_{n})$ be a given $n$-dimensional random vector on a $G$-expectation space $(\Omega, \mathcal{H}, \hat{\mathbb{E}}) .$ We define the functional $\mathbb{F}_{X}$ on the space of Lipschitz functions $C_{L i p}(\mathbb{R}^{n})$ by putting$$
\mathbb{F}_{X}[\varphi]:=\hat{\mathbb{E}}[\varphi(X)], ~ \varphi \in C_{L i p}(\mathbb{R}^{n}).
$$
The triple $(\mathbb{R}^{n}, C_{l . L i p}(\mathbb{R}^{n}), \mathbb{F}_{X})$ forms a nonlinear expectation space, and $\mathbb{F}_{X}$ is called the distribution of $X$ under $\hat{\mathbb{E}}$.
\end{definition}
We also shall introduce the space

$\begin{array}{lll}
M_{G}^{p,0}(0, T )&=&\{\eta_{s}(\omega)=\displaystyle{\sum}_{i=0}^{N-1} \xi_{i}(\omega) I_{(s_{i}, s_{i+1}]}(s): N\geq1,~s_{0}< \cdots< s_{N}\\
& & \qquad\mbox{ partition of } [0,T],\ \xi_{i} \in L_{G}^p(\Omega_{s_{i}}), 0\le i\le N-1\}.
\end{array}$

\noindent By 
$M_{G}^{p}(0, T )$ and $H_{G}^{p}(0, T )$ we denote the completion of $M_{G}^{p,\thinspace 0}(0, T )$ under the norm $\displaystyle\| \cdot \|_{M_{G}^{p}} : = (\hat{\mathbb{E}}[\int_{0}^{T} |\cdot |^{p} d s])^{\frac{1}{p}}$ and $\displaystyle\| \cdot \|_{H_{G}^{p}} : = (\hat{\mathbb{E}}[(\int_{0}^{T} |\cdot |^{2} d s)^{\frac{p}{2}}])^{\frac{1}{p}}$, respectively.

Define {$S_{G}^{0}(0, T):=\{\eta_{s}:=h(s, B_{s_{1} \wedge s}, \ldots, B_{s_{n} \wedge s}): s_{1}, \ldots, s_{n} \in[0, T],\thinspace  h \in C_{b, \text { Lip }}(\mathbb{R}^{n+1})\}$.}

\noindent For $p \geq 1$, we denote by $S_{G}^{p}(0, T)$ the
completion of $S_{G}^{0}(0, T)$ under the norm $\displaystyle\|\eta\|_{S_{G}^{p}}\!\!:=\!\!(\hat{\mathbb{E}}[\!\sup _{s \in[0, T]}\!\!|\eta|^{p}\!])^{\frac{1}{p}}\!\!,$ $\eta\in S_G^0(0,T)$.

Let us now recall the stochastic integration under the $G$-expectation. We  define $\displaystyle\int_{0}^{t} \eta_{s}^n d B_{s}:= \displaystyle{ \sum_{i=0}^{n-1}}\xi_i^n(B_{t_{i+1}}-B_{t_{i}})$, for $\eta_t^n=\displaystyle{ \sum_{i=0}^{N_n-1}}\xi_i^nI_{(t_{i}, t_{i+1}]}(t)\in M_{G}^{2,\thinspace 0}(0, T )$, and for $\eta \in M_{G}^{2}(0, T)$ with $\|\eta^n-\eta\|_{M_{G}^{2}}{\rightarrow}0$ $(n\rightarrow\infty)$, we define  $$\displaystyle\int_{0}^{t} \eta_{s} d B_{s} := L_G^2-\displaystyle{\lim_{n\rightarrow\infty}} \int_{0}^{t} \eta_{s}^n d B_{s},$$ where  $L_G^2$ indicates the convergence in $L^2_G(\Omega)$: $\displaystyle\hat{\mathbb{E}}\Big[\Big|\int_{0}^{t} \eta_{s}^n d B_{s}-\int_{0}^{t} \eta_{s} d B_{s}\Big|^2\Big]{\rightarrow}0$, $n\rightarrow\infty$.

Similarly, we define $\displaystyle\int_{0}^{t} \xi_{s} d\langle B\rangle_{s}$ and $\displaystyle\int_{0}^{t} \xi_{s} d s$ for
 $\xi \in M_{G}^{1}(0, T)$, where $\langle B\rangle$ denotes the cross-variation process of $B$.

Last not least we recall that, given a measurable space  $(\mathbb{X}, \mathscr{X})$ and an $\mathbb{X}$-valued random variable $\xi$ defined on $(\Omega, \mathcal{B}(\mathbb{R}^d), P)$, we denote by $P_{\xi} := P \circ \xi^{-1}$ the law induced by $\xi$ on $(\mathbb{X}, \mathscr{X}) .$
\subsection{SDEs and BSDEs driven by $G$-Brownian motion}
For simplicity, we only consider the one-dimensional case $d=1$, and so also the $G$-Brownian motion is supposed to be one-dimensional. Recall that in this one-dimensional case $G(a)=\frac12 \hat{\mathbb{E}}[aB_1^2]$, and for $\overline{\sigma}^2:=\hat{\mathbb{E}}[B_1^2]$ and $\underline{\sigma}^2:=-\hat{\mathbb{E}}[-B_1^2]$ , we have $G(a)=\frac12\big(\overline{\sigma}^2a^+-\underline{\sigma}^2a^-\big). $ Let us suppose throughout what follows that $\underline{\sigma}^2>0$, i.e., we have $0<\underline{\sigma}^2\le \overline{\sigma}^2<+\infty$. When $\underline{\sigma}^2=\overline{\sigma}^2$, the $G$-expectation is just a linear expection.

We consider the following $G$-SDE: For  given $0 \leq t \leq T<\infty$,
\begin{equation}\label{GSDE}
\left\{\begin{array}{l}
d X_{s}^{t, x}=b(s, X_{s}^{t, x}) d s+ h(s, X_{s}^{t, x}) d\langle B\rangle_{s}+ \sigma(s, X_{s}^{t, x}) d B_{s}, \quad s \in[t, T], \\
X_{t}^{t, x}=x,
\end{array}\right.
\end{equation}
where $x \in \mathbb{R}$, and $b, ~h,~ \sigma: [0, T] \times \Omega \times \mathbb{R} \rightarrow \mathbb{R}$ are given  functions satisfying the following assumptions:

(H1) For some $p\geq2$ it holds $b(\cdot,x),~h(\cdot,x),~\sigma(\cdot,x) \in M_G^p(0,T)$, for all $x\in\mathbb{R}$;

(H2) There exists a constant $L>0$ such that for all $x, x' \in \mathbb{R}$, $t\in[0,T]$,
$$
|b(t,x)-b(t,x')|+|h(t,x)-h(t,x')|+|\sigma(t,x)-\sigma(t,x')| \leq L|x-x'|.
$$
For simplicity, $X_{s}^{0, x}$ will be denoted by $X_{s}^{x}$, for $s\in[0,T],\, x \in \mathbb{R}$. We have the following estimates for $G$-SDE (\ref{GSDE}) which can be found in \cite{11}.
\begin{lemma}\label{lem2.9}
Assume that the conditions $(\rm H 1)$ and $(\rm H 2)$ hold. Then G-SDE (\ref{GSDE}) has a unique solution $(X_{s}^{t, x})_{s \in[t, T]} \in M_{G}^{p}(t, T )$. Moreover, there exists a constant $C\in\mathbb{R}$ depending on $p, T, L$ and $G$ such that, for all $x,~ y \in \mathbb{R},~ t, ~t' \in[0, T]$, we have\\
$$\displaystyle {\rm{i)}}\ \ \hat{\mathbb{E}}\Big[\sup _{s \in[0, t]}|X_{s}^{x}|^{p}\Big] \leq C(1+|x|^{p});\ \ \  {\rm{ii)}}\ \
\hat{\mathbb{E}}[|X_{t}^{x}-X_{t^{\prime}}^{y}|^{p}] \leq C\left(|x-y|^{p}+\left(1+|x|^{p}\right)\left|t-t^{\prime}\right|^{p / 2}\right).
$$
\end{lemma}
We also consider the following BSDE driven by a $G$-Brownian motion:
\begin{equation}\label{1}
Y_{t}=\xi+\int_{t}^{T} f(s, Y_{s}, Z_{s}) d s-\int_{t}^{T} Z_{s} d B_{s}-(K_{T}-K_{t}),~ 0\leq t\leq T,
\end{equation}
where the coefficient
$
f(t, \omega, y, z):[0, T] \times \Omega_{T} \times \mathbb{R}\times \mathbb{R} \rightarrow \mathbb{R}
$
is supposed to satisfy the following conditions:

(H3) There exists some $\beta>1$ such that, for all $y, ~z, ~f(\cdot, \cdot, y, z) \in M_{G}^{\beta}(0, T)$;

(H4) $|f(t, \omega, y, z)-f(t, \omega, y^{\prime}, z^{\prime})| \leq L(|y-y^{\prime}|+|z-z^{\prime}|)$, $(t, \omega)\in[0,T]\times\Omega$, $y,z,y',z'\in\mathbb{R}$, $~~~~~~~~~~~~~~$for some constant $L>0$.

For simplicity, we denote by $\mathfrak{S}_{G}^{p}(0, T)$ the collection of all processes $(Y, Z, K)$ such that $Y \in S_{G}^{p}(0, T), ~Z \in H_{G}^{p}(0, T)$, and $K$ is a non-increasing $G$-martingale with $K_{0}=0$ and $K_{T} \in L_{G}^{p}(\Omega_{T}).$
\begin{definition}[\!\!\cite{12}]\label{def2.10}\rm
Let $\xi \in L_{G}^{\beta}(\Omega_{T})$ and $f$ satisfy (H3) and (H4) for  $\beta>1$. A triplet of processes $(Y, Z, K)$ is called a solution of  (\ref{1}), if for some $1<p \leq \beta$ the following properties hold:

(a) $(Y, Z, K) \in \mathfrak{S}_{G}^{p}(0, T)$;

(b) $\displaystyle Y_{t}=\xi+\int_{t}^{T} f(s, Y_{s}, Z_{s}) d s-\int_{t}^{T} Z_{s} d B_{s}-(K_{T}-K_{t})$, $t\in[0,T]$.
\end{definition}

\begin{theorem}[\!\!\cite{12}]\label{th2.11}
 Assume that $\xi \in L_{G}^{\beta}(\Omega_{T})$ and $f$ satisfies $(\rm H 3)$ and $(\rm H 4)$ for  $\beta>1$. Then (\ref{1}) has a unique solution $(Y, Z, K) .$
 \end{theorem}
\section{Formulation of the Problem}

We consider as  control state space $U$  a non-empty, closed and convex bounded subset of $\mathbb{R}^d$. A process $u:[0,T]\times\Omega\rightarrow U$ is said to be an admissible control on $[0,T]$, if $u\in M_G^2(0,T;U)$. By $\mathcal{U}\big(=M_G^2(0,T;U)\big)$ we denote the class of all admissible controls $u$. For any $u \in \mathcal{U}$, we consider the following stochastic differential equation
\begin{equation}\label{2}\left\{\begin{array}{l}
d X_{t}^{u}=\sigma(X_{t}^{u}, \hat{\mathbb{E}}[\varphi_{1}(X_{t}^{u})], u_{t}) d B_{t}+b(X_{t}^{u} , \hat{\mathbb{E}}[\varphi_{2}(X_{t}^{u})], u_{t}) d t
\\ ~~~~~~~~~
+\beta(X_{t}^{u}, \hat{\mathbb{E}}[\varphi_{3}(X_{t}^{u})], u_{t}) d \langle B \rangle_{t},~ t \in[0, T], \\
X_{0}^{u}=x\in \mathbb{R},
\end{array}\right.
\end{equation}
where
$b,~\beta:[0, T] \times \mathbb{R} \times \mathbb{R} \times U \longrightarrow \mathbb{R},$
$\sigma:[0, T] \times \mathbb{R} \times \mathbb{R} \times U \longrightarrow \mathbb{R},$ and $\varphi_1,~\varphi_2,~\varphi_3: \mathbb{R} \longrightarrow \mathbb{R}.
$

The associated cost functional is given by
\begin{equation}\label{eq3.2}
J(u):=\hat{\mathbb{E}}[\Phi(X_{T}^{u},\hat{\mathbb{E}}[\varphi_{4}(X_{T}^{u})])+\int_{0}^{T} l(t, X_{t}^{u}, \hat{\mathbb{E}}[\varphi_{5}(X_{t}^{u})], u_{t}) d t],
\end{equation}
where
$\Phi: \mathbb{R} \times \mathbb{R} \longrightarrow \mathbb{R},$
$l:[0, T] \times \mathbb{R} \times \mathbb{R} \times U \longrightarrow \mathbb{R},$ and
$\varphi_4,~\varphi_5: \mathbb{R} \longrightarrow \mathbb{R} .$

The following assumptions will be in force throughout this paper.

\noindent(\textbf{A.1}) The functions $\varphi_i,~ i=1,2,3,4,5$,  are continuously differentiable, $\Phi$ and $l$ are continuously \\$~~~~~~~~~$differentiable w.r.t. $(x, y) $, and $b,~ \sigma, ~\beta$ are continuously differentiable w.r.t. $(x, y, v)$.

\noindent(\textbf{A.2}) All the derivatives  in (\textbf{A.1}) are Lipschitz continuous and bounded.

For given $u(\cdot) \in \mathcal{U}$, $X^u$ is called a solution of the above mean-field $G$-SDE  if $X^{u} \in M_{G}^{2}(0, T ; \mathbb{R}^{n}) $ satisfies (\ref{2}).
 Under the
above assumptions, due to Lemma 2.9, SDE (\ref{2}) has a unique solution.
\begin{lemma}[\!\!\cite{19}]\label{th3.1}
(Existence and uniqueness of the solution) If $(\textbf{\rm A.1})$ and $(\textbf{\rm A.2})$ are satisfied, then (\ref{2}) has a unique solution $X^u$, for all $u\in\mathcal{U}$.
\end{lemma}

The optimal control problem consists in minimizing the functional $J(\cdot)$ over $\mathcal{U}$. An admissible control that minimizes $J$ is called optimal.

Our main objective is to characterise the optimal control with the help of Pontryagin's stochastic maximum principle. For this the study of the derivative under the sublinear $G$-expectations is crucial. This is the subject of the following section.

\section{Derivative of a
function of a law under $G$-expectation}

According to Section 2, $(\Omega,\mathcal{H},\hat{\mathbb{E}})$ is a sublinear expectation space, where we restrict now to $\Omega=\Omega_T=C([0,T];\mathbb{R})$. Recall that, due to Theorem \ref{th2.7}, $\mathcal{P}=\{P \text{~a probability on~} (\Omega, \mathcal{B}(\Omega)): E_{P}[X] \leq \hat{\mathbb{E}}[X],
\text{~for~} X \in L_G^1(\Omega)\}$
is a non empty convex, weakly compact subset of $\mathcal{P}(\mathbb{R})$ endowed with the topoplogy of weak convergence. Moreover,
$$
\hat{\mathbb{E}}[\xi]=\sup _{P \in \mathcal{P}} E_{P}[\xi],  \text { for all } \xi \in L_G^1(\Omega),
$$
where the supremum is in fact a maximum: For all $\xi\in L_G^1(\Omega)$, there exists $P\in\mathcal{P}$ such that $\hat{\mathbb{E}}[\xi]=E_{P}[\xi]$ (see \cite{19}). Consequently, the set
$$\mathcal{P}_{\{\xi\}}:=\{P\in \mathcal{P}:\hat{\mathbb{E}}[\xi]=E_{P}[\xi]\}$$
is nonempty.

Let $\xi, \eta\in L_G^1(\Omega)$ and put $F(\lambda):=\hat{\mathbb{E}}[\xi+\lambda\eta], ~\lambda\in \mathbb{R}.$ Now we study the differentiability of $F$. From the definition of the $G$-expectation $\hat{\mathbb{E}}$ we know that   $F$ is convex. Indeed, for all $\lambda, \lambda' \in\mathbb{R}$ and $\rho \in(0,1)$,
$$
\begin{aligned}
&F\big(\rho \lambda+(1-\rho) \lambda'\big)=\hat{\mathbb{E}}\big[\xi+\big(\rho \lambda+(1-\rho) \lambda'\big)\eta\big]
=\hat{\mathbb{E}}\big[\rho\big(\xi+\lambda\eta\big)+(1-\rho)\big(\xi+\lambda'\eta\big)\big]\\
\leq&\rho \hat{\mathbb{E}}\big[\xi+\lambda\eta\big]+(1-\rho) \hat{\mathbb{E}}\big[\xi+\lambda'\eta\big]=\rho F(\lambda)+(1-\rho) F(\lambda').
\end{aligned}
$$
Consequently, for all $\lambda\in \mathbb{R},$ there exists the right-derivative of $F$ at $\lambda$
$$F_{+}^{'}(\lambda)=\lim_{0<\varepsilon\downarrow0}\frac{F(\lambda+\varepsilon)-F(\lambda)}{\varepsilon}$$
and also the corresponding left-derivative
$$F_{-}^{'}(\lambda)=\lim_{0>\varepsilon\uparrow0}\frac{F(\lambda+\varepsilon)-F(\lambda)}{\varepsilon},$$
and, for all $\lambda<\lambda'$, we have
$$F_{-}^{'}(\lambda)\leq F_{+}^{'}(\lambda)\leq\frac{F(\lambda')-F(\lambda)}{\lambda'-\lambda}\leq F_{-}^{'}(\lambda').$$
Let us compute $F_{+}^{'}(0)$ with avoiding the $G$-martingale representation (Recall the $G$-martingale representation from Theorem  \ref{th2.11}, obtained for $f=0$.  Let us also mention that the derivative with use of the $G$-martingale representation as essential tool was discussed in \cite{14}). To this end, we first give the following lemma.
\begin{lemma}\label{lem4.1}
Let $\xi, \eta\in L_G^1(\Omega)$ and $0<\varepsilon_l \downarrow0$ $(l\rightarrow\infty)$, and let $P_{l}\in\mathcal{P}_{\{\xi+\varepsilon_l \eta\}}$, $l\geq1.$ Then we have

$\rm i)$ There exists a subsequence of $(P_l)$, denoted by $(P_{l_k})$, and $P\in \mathcal{P}$, such that  $P_{l_k}\rightharpoonup P$, as $l_k\rightarrow \infty$ (weak convergence of probability measures);

$\rm ii)$ If $P_{l}\rightharpoonup P$, as $\varepsilon_l\downarrow0$ $(l\rightarrow\infty)$, for some $P\in \mathcal{P}$, then $P\in\mathcal{P}_{\{\xi\}}$.
\end{lemma}
\noindent\emph{Proof.} i) From the weak compactness of $\mathcal{P}$ we get i).

ii) Assume that $P_{l}\rightharpoonup P$, as $l\rightarrow\infty$, for some $P\in\mathcal{P}$. Note that the functions in $L_{ip}(\Omega)$ are bounded and uniformly continuous. Thus, for all $\theta\in L_{ip}(\Omega)$, $E_{P_l}[\theta]\rightarrow E_{P}[\theta]$, as $l\rightarrow\infty$. Given any $\delta>0$, let $\theta\in L_{ip}(\Omega)$ be such that  $\hat{\mathbb{E}}[|\theta-\xi|]\leq\delta.$ Then, as $E_{P_l}[\theta]\rightarrow E_{P}[\theta]$, $l\rightarrow\infty$, we have
\begin{equation*}
\begin{aligned}
\big|E_{P_l}[\xi]-E_{P}[\xi]\big|&\leq\big|E_{P_l}[\theta]-E_{P}[\theta]\big|+\big|E_{P_l}[\xi-\theta]\big|+\big|E_{P}[\xi-\theta]\big|\\
&\leq2\delta+\big|E_{P_l}[\theta]-E_{P}[\theta]\big|\rightarrow2\delta, \text{~as~} l\rightarrow\infty.
\end{aligned}
\end{equation*}
From the arbitrariness of $\delta>0$, it follows that $E_{P_l}[\xi]\rightarrow E_{P}[\xi]$, as $l\rightarrow\infty$. But, as $P_l\in\mathcal{P}_{\{\xi+\varepsilon_l\eta\}},\, l\ge 1,$
\begin{equation*}
\begin{aligned}
\big|\hat{\mathbb{E}}[\xi+\varepsilon_l \eta]-              E_{P_l}[\xi]\big|=\big|E_{P_l}[\xi+\varepsilon_l \eta]-E_{P_l}[\xi]\big|\leq \varepsilon_l E_{P_l}[|\eta|]
\leq \varepsilon_l\hat{\mathbb{E}}[|\eta|]\rightarrow0,
\text{~as~}l\rightarrow\infty,
\end{aligned}
\end{equation*}
and also
\begin{equation*}
\begin{aligned}
\big|\hat{\mathbb{E}}[\xi+\varepsilon_l \eta]-\hat{\mathbb{E}}[\xi]\big|\leq\varepsilon_l \hat{\mathbb{E}}[|\eta|]\rightarrow0, \text{~as~}l\rightarrow\infty,
\end{aligned}
\end{equation*}
it follows that $\hat{\mathbb{E}}[\xi]=E_P[\xi]$, i.e., $P\in\mathcal{P}_{\{\xi\}}.$\hfill$\square$

We recall that the set $\mathcal{P}$ endowed with the weak convergence of probability measures is a compact metrisable space. Let $d(\cdot,\cdot)$ be a metric on $\mathcal{P}$ which is compatible with the weak convergence, e.g., we can choose the L\'evy-Prokhorov metric (see Theorem 11.3-3, \cite{25}):

\centerline{$\displaystyle d(P,Q):=\sup\big\{\int_{\Omega}fd(P-Q),\, |f|_{BL}\le 1 \big\},$}

\noindent where $\displaystyle |f|_{BL}=\sup_{\omega\in\Omega}|f(\omega)|+\sup_{\omega\not=\omega'}\frac{|f(\omega)-f(\omega')|}{\ \ \  |\omega-\omega'|_{C([0,T])}}$.

Observe also that, as $(\mathcal{P},d)$ is a compact metric space, it is, in particular, also separable.

For $A, B\subset \mathcal{P}$, we put
$$d(P,B):=\mbox{dist}_B(P)=\inf \{d(P, Q): Q \in B\}, \text{~for~}  P \in \mathcal{P}, \text{~and~} \Gamma(A, B):=\sup_{P\in A} d(P,B).$$
Note that $\Gamma(A,B)$ is the maximal distance from $B$ of the probabilities in $A$. In particular, $\Gamma(A,B)=0,$ if $A\subset B$. Of course, $\Gamma(\cdot,\cdot)$ is not symmetric, its symmetrisation is just the Hausdorff distance $d_{H}(A, B)=\max \{\Gamma(A, B), \Gamma(B, A)\},~ A, B \subset \mathcal{P}$.
\begin{prop}\label{prop4.2} We have
$\Gamma\left(\mathcal{P}_{\{\xi+\varepsilon\eta\}}, \mathcal{P}_{\{\xi\}}\right) \rightarrow 0$, as $\varepsilon\downarrow 0$.
\end{prop}
\noindent\emph{Proof.} Let $0<\varepsilon_{l} \downarrow 0$ and $P_{l} \in \mathcal{P}_{\{\xi+\varepsilon_l\eta\}}$ be such that
$$\Gamma\left(\mathcal{P}_{\{\xi+\varepsilon_l\eta\}}, \mathcal{P}_{\{\xi\}}\right)-\frac{1}{l} \leq d\left(P_{l}, \mathcal{P}_{\{\xi\}}\right), ~l \geq 1.$$
Due to Lemma \ref{lem4.1}, for all subsequence $(P_{l_{k}})_{k \geq 1} \subset(P_{l})_{l \geq 1}$, there exists some sub-subsequence $(P_{l_{k_n}})_{n \geq1}\subset(P_{l_{k}})_{k \geq 1}$ and some $P \in \mathcal{P}_{\{\xi\}}$ such that $P_{l_{k_{n}}} \rightharpoonup P$, as $n \rightarrow \infty$. Then, $$\Gamma\big(\mathcal{P}_{\{\xi+\varepsilon_{l_{k_n}}\eta\}}, \mathcal{P}_{\{\xi\}}\big) \leq d\left(P_{l_{k_{n}}}, \mathcal{P}_{\{\xi\}}\right)+\frac{1}{l_{k_{n}}} \leq d\left(P_{l_{k_{n}}}, P\right)+\frac{1}{l_{k_{n}}} \rightarrow 0, \text{~as~} n \rightarrow \infty.$$
This implies
$$\Gamma\left(\mathcal{P}_{\{\xi+\varepsilon_l\eta\}}, \mathcal{P}_{\{\xi\}}\right) \rightarrow 0~~(l \rightarrow \infty),$$
for any $0<\varepsilon_{l} \downarrow 0$, and, consequently,
$\Gamma\left(\mathcal{P}_{\{\xi+\varepsilon\eta\}}, \mathcal{P}_{\{\xi\}}\right) \rightarrow 0, \text{~as~} \varepsilon\downarrow 0.$\hfill$\square$
\begin{remark}\rm Lemma \ref{lem4.1} can also be regarded as a consequence of Proposition \ref{prop4.2}.
Indeed, for any $0<\varepsilon_{l} \downarrow 0$, let $P_{l} \in \mathcal{P}_{\{\xi+\varepsilon_l\eta\}},~ l \geq 1$. Then,
$$d\left(P_{l}, \mathcal{P}_{\{\xi\}}\right) \leq \Gamma\left(\mathcal{P}_{\{\xi+\varepsilon_l\eta\}}, \mathcal{P}_{\{\xi\}}\right) \rightarrow 0,\text{~as~}0<\varepsilon_{l} \downarrow 0.$$
Let $Q_{l} \in \mathcal{P}_{\{\xi\}}$ be such that $d\left(P_{l}, Q_{l}\right) \leq d\left(P_{l}, \mathcal{P}_{\{\xi\}}\right)+\frac{1}{l},~ l \geq 1$.
As $\mathcal{P}$ is weakly compact, there is a subsequence $(Q_{l_{k}})_{k \geq 1}\subset(Q_{l})_{l \geq 1}$ and some $Q \in \mathcal{P}$ such that $Q_{l_{k}} \rightharpoonup Q$. Consequently, due to the corresponding argument in the proof of Lemma \ref{lem4.1}, $\hat{\mathbb{E}}[\xi]=E_{Q_{l_{k}}}[\xi] \rightarrow E_{Q}[\xi]$, i.e., also $Q \in \mathcal{P}_{\{\xi\}}$. Finally, from
$
d(P_{l_{k}}, Q) \leq d(P_{l_{k}}, Q_{l_{k}})+d(Q_{l_{k}}, Q) \leq \Gamma\big(\mathcal{P}_{\{\xi+\varepsilon_{l_k}\eta\}}, \mathcal{P}_{\{\xi\}}\big)+\frac{1}{l_k}+d ( Q_{l_{k}}, Q)\rightarrow 0, ~k \rightarrow \infty,
$
 we see that $P_{l_{k}} \rightharpoonup Q \in \mathcal{P}_{\{\xi\}}$, as $k \rightarrow \infty$.
\end{remark}

\medskip

Let us now come to the computation of the right-derivative $F'_+(0)$ of $F(\lambda)=\hat{\mathbb{E}}[\xi+\lambda\eta]$ ($\xi,\, \eta\in  L^1_G(\Omega)$) at $\lambda=0$. For this we let  $0<\varepsilon_l \downarrow0$, $P_{l}\in\mathcal{P}_{\{\xi+\varepsilon_l \eta\}}$ and $P\in\mathcal{P}_{\{\xi\}}$ be such that $P_l\rightharpoonup P$ (Due to Lemma \ref{lem4.1} this choice is possible). In analogy to the fact that $P_l\rightharpoonup P$ implies $E_{P_{l}}[\xi]\rightarrow E_P[\xi]$, we get that, for any $\zeta\in L_G^1(\Omega)$, $E_{P_{l}}[\zeta]\rightarrow E_P[\zeta]$, as
$l\rightarrow\infty$, and so $E_{P_{l}}[\eta]\rightarrow E_P[\eta]$, as $l\rightarrow\infty$. Then, as $P_{l}\in\mathcal{P}_{\{\xi+\varepsilon_l \eta\}}$,
$$
F_{+}^{'}(0)=\lim_{\varepsilon\downarrow0}\frac{\hat{\mathbb{E}}[\xi+\varepsilon \eta]-\hat{\mathbb{E}}[\xi]}{\varepsilon}\leq\lim_{l\rightarrow\infty}\frac{E_{P_{l}}[\xi+\varepsilon_l \eta]-E_{P_{l}}[\xi]}{\varepsilon_l}=E_{P}[\eta],
$$
i.e., $F_{+}^{'}(0)\leq E_{P}[\eta]$.
On the other hand, for all $Q\in\mathcal{P}_{\{\xi\}}$,
$$
F_{+}^{'}(0)=\lim_{\varepsilon\downarrow0}\frac{\hat{\mathbb{E}}[\xi+\varepsilon \eta]-\hat{\mathbb{E}}[\xi]}{\varepsilon}\geq\lim_{\varepsilon\downarrow0}\frac{E_{Q}[\xi+\varepsilon \eta]-E_{Q}[\xi]}{\varepsilon}=E_{Q}[\eta].
$$
Consequently, we get the following lemma.
\begin{lemma}\label{lem4.4}
For $\xi, \eta\in L_G^1(\Omega)$ and $F(\lambda):=\hat{\mathbb{E}}[\xi+\lambda\eta]$, we have $$F_{+}^{'}(0)=\displaystyle{\max_{P\in\mathcal{P}_{\{\xi\}}}}E_P[\eta]=\hat{\mathbb{E}}_{\{\xi\}}[\eta],$$
where $\hat{\mathbb{E}}_{\{\xi\}}[\eta]:=\displaystyle{\sup_{P\in\mathcal{P}_{\{\xi\}}}}E_{P}[\eta]$, $\eta\in L_G^1(\Omega)$, is a new sublinear expectation, and $\hat{\mathbb{E}}_{\{\xi\}}$ is dominated by $\hat{\mathbb{E}}$, i.e., $\hat{\mathbb{E}}_{\{\xi\}}[\thinspace\cdot\thinspace]\leq\hat{\mathbb{E}}[\thinspace\cdot\thinspace]$.
\end{lemma}
\begin{remark}\rm
 From the above lemma it follows that
$$
F_{-}^{'}(0)=\lim_{0<\varepsilon\downarrow0}\frac{\hat{\mathbb{E}}[\xi-\varepsilon \eta]-\hat{\mathbb{E}}[\xi]}{-\varepsilon}=-\lim_{0<\varepsilon\downarrow0}\frac{\hat{\mathbb{E}}[\xi+\varepsilon (-\eta)]-\hat{\mathbb{E}}[\xi]}{\varepsilon}=-\hat{\mathbb{E}}_{\{\xi\}}[-\eta].
$$
This shows in particular that $F(\lambda)=\hat{\mathbb{E}}[\xi+\lambda\eta]$ is differentiable at $\lambda=0$ if and only if
$\hat{\mathbb{E}}_{\{\xi\}}[\eta]=-\hat{\mathbb{E}}_{\{\xi\}}[-\eta].$

We also observe that, for all $\lambda\in\mathbb{R}$,
$$
F_{+}^{'}(\lambda)=\lim_{0<\varepsilon\downarrow0}\frac{\hat{\mathbb{E}}[(\xi+\lambda\eta)+\varepsilon \eta]-\hat{\mathbb{E}}[\xi+\lambda\eta]}{\varepsilon}=\hat{\mathbb{E}}_{\{\xi+\lambda\eta\}}[\eta],
$$
$$
F_{-}^{'}(\lambda)=-\lim_{0<\varepsilon\downarrow0}\frac{\hat{\mathbb{E}}[(\xi+\lambda\eta)+\varepsilon (-\eta)]-\hat{\mathbb{E}}[\xi+\lambda\eta]}{\varepsilon}=-\hat{\mathbb{E}}_{\{\xi+\lambda\eta\}}[-\eta],
$$
and as $F$ is convex, for all $\lambda<\lambda'$,
$$-\hat{\mathbb{E}}_{\{\xi+\lambda\eta\}}[-\eta]\leq\hat{\mathbb{E}}_{\{\xi+\lambda\eta\}}[\eta]\leq\frac{\hat{\mathbb{E}}[\xi+\lambda'\eta]-\hat{\mathbb{E}}[\xi+\lambda\eta]}{\lambda'-\lambda}\leq-\hat{\mathbb{E}}_{\{\xi+\lambda'\eta\}}[-\eta].$$
\end{remark}
\begin{corollary}\label{coro4.5}
Let $\varphi\in C^1(\mathbb{R})$ have a bounded Lipschitz derivative $\partial\varphi:\mathbb{R}\rightarrow\mathbb{R}$, and let $\xi,\ \eta\in L_G^1(\Omega)$. Then, for $H(\lambda):= \hat{\mathbb{E}}[\varphi(\xi+\lambda\eta)], ~\lambda\in\mathbb{R}$, we have
$$ {\rm i)}\  H^{'}_{+}(0)=\hat{\mathbb{E}}_{\{\varphi(\xi)\}}[\partial\varphi(\xi)\eta];\ \  and\ \ {\rm ii)}\  H^{'}_{-}(0)=-\hat{\mathbb{E}}_{\{\varphi(\xi)\}}[-\partial\varphi(\xi)\eta].$$
\end{corollary}
\noindent\emph{Proof.} Let $\xi':=\varphi(\xi)$ and $\eta':=\partial\varphi(\xi)\eta$. Then, $\xi', \eta'\in L_G^1(\Omega)$, and
\begin{equation*}
\begin{aligned}
\varphi(\xi+\varepsilon\eta)&=\varphi(\xi)+\int_0^1 \partial_{\lambda}[\varphi(\xi+\lambda\varepsilon\eta)]d\lambda=\varphi(\xi)+\int_0^1 \partial\varphi(\xi+\lambda\varepsilon\eta)d\lambda\cdot\varepsilon\eta\\
&=\varphi(\xi)+\varepsilon\big(\partial\varphi(\xi)\eta\big)+R_{\varepsilon},
\end{aligned}
\end{equation*}
where
$\displaystyle
R_{\varepsilon}:=\int_0^1\Big(\partial\varphi(\xi+\lambda\varepsilon\eta)-\partial\varphi(\xi)\Big)d\lambda\cdot\varepsilon\eta,
$
and from the Lipschitz property of $\partial\varphi$ we have\\
$\displaystyle |R_{\varepsilon}|\leq C\varepsilon^2|\eta|^2$\ and\ $\hat{\mathbb{E}}[|R_{\varepsilon}|]\leq C_{\eta}\varepsilon^2,~\varepsilon>0.$
Hence, as $\varphi(\xi+\varepsilon\eta)=\xi'+\varepsilon\eta'+R_{\varepsilon}$,
\begin{equation*}
\begin{aligned}
\bigg|\frac{\hat{\mathbb{E}}[\varphi(\xi+\varepsilon\eta)]-\hat{\mathbb{E}}[\varphi(\xi)]}{\varepsilon}-\hat{\mathbb{E}}_{\{\varphi(\xi)\}}[\partial\varphi(\xi)\eta]\bigg|\leq\bigg|\frac{\hat{\mathbb{E}}[\xi'+\varepsilon\eta']-\hat{\mathbb{E}}[\xi']}{\varepsilon}-\hat{\mathbb{E}}_{\{\xi'\}}[\eta']\bigg|+\hat{\mathbb{E}}[\frac{1}{\varepsilon}|R_{\varepsilon}|]\rightarrow0,
\end{aligned}
\end{equation*}
as $\varepsilon\downarrow0$. This proves i). For ii) we note that, from i) with $-\eta$ instead of $\eta$, it follows that
$$~~~~~~~~~~~~~~~~
H^{'}_{-}(0)=-\lim_{0<\varepsilon\downarrow0}\frac{\hat{\mathbb{E}}[\varphi(\xi+\varepsilon(-\eta))]-\hat{\mathbb{E}}[\varphi(\xi)]}{\varepsilon}=-\hat{\mathbb{E}}_{\{\varphi(\xi)\}}[-\partial\varphi(\xi)\eta].
~~~~~~~~~~~~~~~~~~~\square$$
\begin{remark}\label{remark4.7} Observe that Corollary \ref{coro4.5} can be extended in a straight-forward way to d-dimensional random variables $\xi,\ \eta\in (L_G^1(\Omega))^d$. Indeed, in the proof of the corollary it suffices the derivative $\partial\varphi$ replaced by the gradient $\nabla\varphi$.
 \end{remark}

\section{Pontryagin's SMP for
mean field stochastic control problems under $G$-expectation}
\subsection{Necessary conditions for optimality}
In the next both sections, to simplify the dynamics and the related computations, we put $\beta = 0$, and so SDE (\ref{2}) becomes
\begin{equation}\label{5.1}\left\{ \begin{array}{l}
dx_{t}^{u}=b(x_{t}^{u} , \hat{\mathbb{E}}[\varphi_{1}(x_{t}^{u})], u_{t}) d t +\sigma(x_{t}^{u}, \hat{\mathbb{E}}[\varphi_{2}(x_{t}^{u})], u_{t}) d B_{t}, \ t \in[0, T],\\
x_{0}^{u}=x\in \mathbb{R}^n.
\end{array}\right.
\end{equation}

\noindent The cost functional is still given by  \eqref{eq3.2},
\begin{equation}
J(u):=\hat{\mathbb{E}}[\Phi(x_{T}^{u},\hat{\mathbb{E}}[\varphi_{4}(x_{T}^{u})])+\int_{0}^{T} l(t, x_{t}^{u}, \hat{\mathbb{E}}[\varphi_{5}(x_{t}^{u})], u_{t}) d t],\ \mbox{for any}\ u\in \mathcal{U}.
\end{equation}

We suppose that there exists an optimal control $\hat{u}\in\mathcal{U}$, that is, $J(\hat{u})\leq J(u)$, for all $u\in\mathcal{U}$. Let us denote $\hat{x}_t:=x_t^{\hat{u}}$, $t\in[0,T]$.
\subsubsection{Taylor expansions}
Let $u$ be an arbitrary admissible control process in $\mathcal{U}$. By $x^{\theta}$ we denote the state process defined by SDE (\ref{5.1}) with the control process $u^{\theta}$ defined as convex perturbation of $\hat{u}$:
$$u_t^{\theta}=\hat{u}_t+\theta(u_t -\hat{u}_t),\ t\in[0,T],\ \theta\in[0,1].$$
We put $v_t=u_t -\hat{u}_t$, $t\in[0,T],$ and introduce the following notations:
\begin{equation}\label{eq5.2'}
\begin{array}{llll}
&\hat{b}(t)=b(t,\hat{x}_t,\hat{\mathbb{E}}[\hat{\varphi}_{1}(t)],\hat{u}_t),
& &\hat{b}_x(t)=b_x(t,\hat{x}_t,\hat{\mathbb{E}}[\hat{\varphi}_{1}(t)],\hat{u}_t),\\
&\hat{\sigma}(t)=\sigma(t,\hat{x}_t,\hat{\mathbb{E}}[\hat{\varphi}_{2}(t)],\hat{u}_t), & &\hat{\sigma}_x(t)=\sigma_x(t,\hat{x}_t,\hat{\mathbb{E}}[\hat{\varphi}_{2}(t)],\hat{u}_t),\\
&\hat{\Phi}(t)=\Phi(\hat{x}_t,\hat{\mathbb{E}}[\hat{\varphi}_{4}(t)]), & &\hat{\Phi}_x(t)=\Phi_x(\hat{x}_t,\hat{\mathbb{E}}[\hat{\varphi}_{4}(t)]),\\
&\hat{l}(t)=l(t,\hat{x}_t,\hat{\mathbb{E}}[\hat{\varphi}_{5}(t)],\hat{u}_t), & &\hat{l}_x(t)=l_x(t,\hat{x}_t,\hat{\mathbb{E}}[\hat{\varphi}_{5}(t)],\hat{u}_t),
\end{array}
\end{equation}
with $\hat{\varphi}_i(t)=\varphi_{i}(\hat{x}_t),\ \hat{\varphi}'_i(t)=\varphi'_{i}(\hat{x}_t), i=1,2,4,5,$
and similarly are defined $\hat{b}_y(t)$, $\hat{b}_v(t)$, $\hat{\sigma}_y(t)$, $\hat{\sigma}_v(t)$, $\hat{l}_y(t)$ and $\hat{l}_v(t)$.
Here, for notational convenience, we denote by $b_{x}$, $b_{y}$, $b_{v}$  the derivative of $b$ w.r.t. the state trajectory, the expected value and the control variable, respectively, and similarly for the other functions.

The objective of this section is to determine the directional derivative of the cost functional in terms of the first order Taylor expansion of the state process. We begin with identifying the Taylor expansion.
\begin{lemma}\label{lem5.1}
 Let $z=(z_t)\in M_G^2(0,T;\mathbb{R})$ be the unique solution of the following SDE

\begin{equation}\label{5.3}
\left\{\begin{aligned}
\mathrm{d} z_{t}=&\Big(\hat{b}_{x}(t) z_{t}+\hat{b}_{y}(t) \hat{\mathbb{E}}_{\{\varphi_1(\hat{x}_{t})\}}[\hat{\varphi}_{1}'(t) z_{t}]+\hat{b}_{v}(t) v_{t}\Big) \mathrm{d} t \\
&+\Big(\hat{\sigma}_{x}(t) z_{t}+\hat{\sigma}_{y}(t) \hat{\mathbb{E}}_{\{\varphi_2(\hat{x}_{t})\}}[\hat{\varphi}_{2}'(t) z_{t}]+\hat{\sigma}_{v}(t) v_{t}\Big) \mathrm{d} B_{t},~t \in[0, T], \\
z_{0}=0 & .
\end{aligned}\right.
\end{equation}
Then, it holds that
$$
\lim _{\theta \rightarrow 0}\hat{\mathbb{E}}[\sup_{t\in[0,T]}|\frac{x_{t}^{\theta}-\hat{x}_{t}}{\theta}-z_{t}|^2]=0.
$$
\end{lemma}
\noindent\emph{Proof.} First we observe that, thanks to our assumptions on the coefficients, we have the existence and the uniqueness for SDE \eqref{5.3}. Now, to simplify our computations, but without loss of generality for the method of the proof,  let $b=0$. So SDE (\ref{5.1}) becomes
\begin{equation}\left\{ \begin{array}{l}
dx_{t}^{u}=\sigma(x_{t}^{u}, \hat{\mathbb{E}}[\varphi_{2}(x_{t}^{u})], u_{t}) d B_{t}, \ \ t \in[0, T],\\
x_{0}^{u}=x,
\end{array}\right.
\end{equation}
while SDE (\ref{5.3}) writes
\begin{equation}\label{5.5}
\left\{\begin{aligned}
\mathrm{d} z_{t}=&\Big(\hat{\sigma}_{x}(t) z_{t}+\hat{\sigma}_{y}(t) \hat{\mathbb{E}}_{\{\varphi_2(\hat{x}_{t})\}}[\hat{\varphi}_{2}'(t) z_{t}]+\hat{\sigma}_{v}(t) v_{t}\Big) \mathrm{d} B_{t},\ \  t \in[0, T], \\
z_{0}=0 & .
\end{aligned}\right.
\end{equation}
Putting $\hat{\rho}_{t}:=\hat{\mathbb{E}}[\varphi_{2}(\hat{x}_{t})]$, $\rho_{t}^{\theta}:=\hat{\mathbb{E}}[\varphi_{2}(x_{t}^{\theta})]$, $u_{t}^{\theta}:=\hat{u}_{t}+\theta v_{t}$ and $\Theta^{\lambda}_t:=(\hat{x}_t+\lambda(x^{\theta}_t-\hat{x}_t),\hat{\rho}_t+\lambda(\rho^{\theta}_t-\hat{\rho}_t),\hat{u}_t+\lambda\theta v_{t})$, we have
\begin{equation}\label{5.6}
\begin{aligned}
&\frac{1}{\theta}(x_{t}^{\theta}-\hat{x}_{t})=\frac{1}{\theta} \int_{0}^{t}\big(\sigma(x_{s}^{\theta}, \rho_{s}^{\theta}, u_{s}^{\theta})-\sigma(\hat{x}_{s}, \hat{\rho}_{s}, \hat{u}_{s})\big) d B_{s}
=\frac{1}{\theta} \int_{0}^{t} \int_{0}^{1} \partial_{\lambda}\big[\sigma(\Theta^{\lambda}_s)] d \lambda d B_{s}\\
=&\frac{1}{\theta} \int_{0}^{t} \int_{0}^{1}\Big\{\sigma_{x}(\Theta^{\lambda}_s)(x^{\theta}_s-\hat{x}_s)+\sigma_{y}(\Theta^{\lambda}_s)(\rho^{\theta}_s-\hat{\rho}_s)+\sigma_{v}(\Theta^{\lambda}_s)\theta v_{s}\Big\}  d \lambda d B_{s}\\
=&\int_{0}^{t}\Big\{\hat{\sigma}_{x}(s) \frac{1}{\theta}(x_{s}^{\theta}-\hat{x}_{s})+\hat{\sigma}_{y}(s) \frac{1}{\theta}(\rho_{s}^{\theta}-\hat{\rho}_{s})+\hat{\sigma}_{v}(s) v_{s}\Big\} d B_{s}+R_{t}^{\theta},\\
\end{aligned}
\end{equation}
where
\begin{equation}\label{5.7}
\begin{aligned}
R_{t}^{\theta}:=\int_{0}^{t}\! \! \int_{0}^{1}\!\!\!\Big\{\!\big(\sigma_{x}(\Theta^{\lambda}_{s})\!-\!\hat{\sigma}_{x}(s)\big) \frac{x_{s}^{\theta}\!-\!\hat{x}_{s}}{\theta}+\big(\sigma_{y}(\Theta^{\lambda}_{s})\!-\!\hat{\sigma}_{y}(s)\big) \frac{\rho_{s}^{\theta}\!-\!\hat{\rho}_{s}}{\theta}+\big(\sigma_{v}(\Theta^{\lambda}_{s})\!-\!\hat{\sigma}_{v}(s)\big) v_{s}\Big\} d \lambda d B_{s}.
\end{aligned}
\end{equation}
We put $\nabla\sigma:=(\sigma_{x}, \sigma_{y}, \sigma_{v})$ and $\nabla \hat{\sigma}(s):=\big(\hat{\sigma}_{x}(s), \hat{\sigma}_{y}(s), \hat{\sigma}_{z}(s)\big)$. Then, thanks to Assumption (A.2),
$$
\big|\nabla \sigma(\Theta^{\lambda}_{s})-\nabla \hat{\sigma}(s)\big| \leq C\big(|x_{s}^{\theta}-\hat{x}_{s}|+|\rho_{s}^{\theta}-\hat{\rho}_{s}|+\theta|v_{s}|\big).
$$
Since $\displaystyle{x_{t}^{\theta}-\hat{x}_{t}=\int_{0}^{t}\big(\sigma(x_{s}^{\theta}, \rho_{s}^{\theta}, u_{s}^{\theta})-\sigma(\hat{x}_{s}, \hat{\rho}_{s}, \hat{u}_{s})\big) d B_{s}}$,  for all $p\geq2$, we have
\begin{equation*}
\begin{aligned}
\hat{\mathbb{E}}\big[\sup_{0\leq s\leq t} |x_{s}^{\theta}-\hat{x}_{s}|^{p}\big] &\leq C_{p} \hat{\mathbb{E}}\Big[\Big(\int_{0}^{t}\big(|x_{s}^{\theta}-\hat{x}_{s}|+|\rho_{s}^{\theta}-\hat{\rho}_{s}|+\theta|v_{s}|\big)^{2} d s\Big)^{\frac{p}{2}}\Big]\\
&\leq C_{p} \hat{\mathbb{E}}\Big[\int_{0}^{t}\Big(|x_{s}^{\theta}-\hat{x}_{s}|^p+\big(\hat{\mathbb{E}}\big[|x_{s}^{\theta}-\hat{x}_{s}|\big]\big)^p+\theta^p|v_{s}|^p\Big) d s\Big]\\
&\leq C_{p}\Big(\theta^{p}+\int_{0}^{t} \hat{\mathbb{E}}\big[|x_{s}^{\theta}-\hat{x}_{s}|^{p}\big] d s\Big),~ t \in[0, T]\\
\end{aligned}
\end{equation*}
(Recall that the control state space $\mathcal{U}$ is bounded).

\noindent Then, by Gronwall's Lemma we have
\begin{equation}\label{5.8}
\hat{\mathbb{E}}\big[\sup_{0\leq s\leq T}|x_{s}^{\theta}-\hat{x}_{s}|^{p}\big] \leq C_p \theta^{p},~ \theta>0.
\end{equation}
Hence, as $\varphi_2$ is Lipschitz, also for $\rho^{\theta}-\hat{\rho}$ we have
\begin{equation}\label{5.9}
\sup_{0\leq s\leq T}|\rho_{s}^{\theta}-\hat{\rho}_{s}|^{p} \leq C_p \theta^{p}, ~\theta>0.
\end{equation}
On the other hand, by standard estimates we have $\displaystyle\hat{\mathbb{E}}\Big[\sup_{0\leq s\leq T} |z_{s}|^{p}\Big] \leq C_{p},~ p \geq 2.$

\noindent From (\ref{5.7}), (\ref{5.8}) and (\ref{5.9}), for some $p \geq 2$,
\begin{equation}\label{ER}
\begin{aligned}
\hat{\mathbb{E}}\big[\sup_{t\in[0,T]} |R_t^{\theta}|^{p}\big] &\leq C_{p} \hat{\mathbb{E}}\Big[\Big(\int_{0}^{T}\Big(\big(|x_{s}^{\theta}-\hat{x}_{s}|+|\rho_{s}^{\theta}-\hat{\rho}_{s}|+\theta|v_{s}|\big)^{2}\frac{1}{\theta}\Big)^2 d s\Big)^{\frac{p}{2}}\Big]\\
&\leq C_{p} \frac{1}{\theta^p} \Big(\hat{\mathbb{E}}\big[\sup_{t\in[0,T]}|x_{s}^{\theta}-\hat{x}_{s}|^{2p}\big]+\sup_{t\in[0,T]}|\rho_{s}^{\theta}-\hat{\rho}_{s}|^{2p}+\theta^{2p}\Big)\leq C_p \theta^{p}, ~\theta>0.
\end{aligned}
\end{equation}
We put ${y_{t}^{\theta}:=\frac{x_{t}^{\theta}-\hat{x}_{t}}{\theta}-z_{t}}$, $t \in[0, T]$, and we have
\begin{equation*}
\begin{aligned}
&\varphi_{2}(x_{s}^{\theta})-\varphi_{2}(\hat{x}_{s})=\Big(\hskip -0.1cm\int_{0}^{1} \hskip -0.1cm\varphi'_{2}\big(\hat{x}_{s}+\lambda(x_{s}^{\theta}-\hat{x}_{s})\big) d \lambda\Big)(x_{s}^{\theta}-\hat{x}_{s})=\varphi'_{2}(\hat{x}_{s})(x_{s}^{\theta}-\hat{x}_{s})+\tilde{R}_{s}^{\theta}\\
=&\ \theta(\varphi'_{2}(\hat{x}_{s})z_{s})+\varphi'_{2}(\hat{x}_{s})(x_{s}^{\theta}-\hat{x}_{s}-\theta z_{s})+\tilde{R}_{s}^{\theta}
=\theta(\varphi'_{2}(\hat{x}_{s})z_{s})+\theta \varphi'_{2}(\hat{x}_{s})y_{s}^{\theta}+\tilde{R}_{s}^{\theta},
\end{aligned}
\end{equation*}
where $\displaystyle\tilde{R}_{s}^{\theta}:=\Big(\hskip -0.1cm \int_{0}^{1}\hskip -0.2cm\big(\varphi'_{2}\big(\hat{x}_{s}+\lambda( x_{s}^{\theta}-\hat{x}_{s})\big)-\varphi'_{2}( \hat{x}_{s}) \big) d \lambda\Big)(x_{s}^{\theta}-\hat{x}_{s})$, and $|\tilde{R}_{s}^{\theta}| \leq C|x_{s}^{\theta}-\hat{x}_{s}|^{2}$, i.e., \begin{equation}\label{R}
\varphi_{2}(x_{s}^{\theta})=\varphi_{2}(\hat{x}_{s})+\theta\big(\varphi'_{2}(\hat{x}_{s}) z_{s}\big)+\theta \varphi'_{2}(\hat{x}_{s}) y_{s}^{\theta}+\tilde{R}_{s}^{\theta}.
\end{equation}
Furthermore, we get
\begin{equation}\label{5.10}
\begin{aligned}
\frac{\rho_{s}^{\theta}-\hat{\rho}_{s}}{\theta}&=\frac{1}{\theta}\Big(\hat{\mathbb{E}}\big[\varphi_{2}(x_{s}^{\theta})\big]-\hat{\mathbb{E}}\big[\varphi_{2}(\hat{x}_{s})\big] \Big)\\
&=\frac{1}{\theta}\Big(\hat{\mathbb{E}}\big[\varphi_{2}(\hat{x}_{s})+\theta\big(\varphi'_{2}(\hat{x}_{s}) z_{s}\big)\big]-\hat{\mathbb{E}}\big[\varphi_{2}(\hat{x}_{s})\big]\Big)+\hat{R}_s^{\theta}(y_s^{\theta}),
\end{aligned}
\end{equation}
where, thanks to (\ref{R}),
$$
\hat{R}_{s}^{\theta}(\eta):=\frac{1}{\theta}\Big(\hat{\mathbb{E}}\big[\varphi_{2}(\hat{x}_{s})+\theta\big(\varphi'_{2}(\hat{x}_{s}) z_{s}\big)+\theta \varphi'_{2}(\hat{x}_{s}) \eta+\tilde{R}_{s}^{\theta}\big]-\hat{\mathbb{E}}\big[\varphi_{2}(\hat{x}_{s})+\theta\big(\varphi'_{2}(\hat{x}_{s}) z_{s}\big)\big]\Big),~
\eta \in L_{G}^{2}(\Omega_{T}).
$$
Notice that,  due to \eqref{5.8},
\begin{equation}\label{RR}{|\hat{R}_{s}^{\theta}(\eta)| \leq C \hat{\mathbb{E}}[|\eta|]+C \hat{\mathbb{E}}\Big[\frac{|x_{s}^{\theta}-\hat{x}_{s}|^{2}}{\theta}\Big]\leq C\big(\theta+\hat{\mathbb{E}}[|\eta|]\big),~ \theta>0.}\end{equation}

Now, we define ${F_{s}(\theta):=\hat{\mathbb{E}}\big[\varphi_{2}(\hat{x}_{s})+\theta\big(\varphi'_{2}(\hat{x}_{s}) z_{s}\big)\big], ~\theta \geq 0}$. Then, as $F_{s}: \mathbb{R}_{+}  \rightarrow \mathbb{R}$ is convex, we obtain
$$0 \leq H_{s}(\theta):=\frac{F_{s}(\theta)-F_{s}(0)}{\theta}-F'_{s, +}(0) \leq F'_{s,+}(\theta)-F'_{s, +}(0) \downarrow 0, \text{~as~} \theta \downarrow 0,$$
where, due to Lemma \ref{lem4.4}, the right-derivative of $F_s(\lambda)$ at $\lambda=0$ satisfies
$$F'_{s,+}(0)=\hat{\mathbb{E}}_{\{\varphi_{2}(\hat{x}_{s})\}}\big[\varphi'_{2}(\hat{x}_{s}) z_{s}\big].$$
Observe that, thanks to (\ref{5.10}) and the above definition of $F_s(\theta)$ and $H_s(\theta)$,
\begin{equation}\label{rho}
\frac{\rho_{s}^{\theta}-\hat{\rho}_{s}}{\theta}=H_{s}(\theta)+ F'_{s, +}(0)+\hat{R}_{s}^{\theta}(y_{s}^{\theta}), ~ s \in[0, T].
\end{equation}
Then, recalling that $y_{t}^{\theta}:=\frac{x_{t}^{\theta}-\hat{x}_{t}}{\theta}-z_{t}$, from (\ref{5.5}) and (\ref{5.6}), we get
\begin{equation*}
\left\{\begin{aligned}
&d y_{t}^{\theta}=d\Big(\frac{1}{\theta}(x_{t}^{\theta}-\hat{x}_{t})-z_{t}\Big)=\Big(\hat{\sigma}_{x}(t) y_{t}^{\theta}+\hat{\sigma}_{y}(t)\Big(\frac{1}{\theta}(\rho_{t}^{\theta}-\hat{\rho}_{t})-\hat{\mathbb{E}}_{\{\varphi_{2}(\hat{x}_{t})\}}\big[\hat{\varphi}'_{2}(t) z_{t}\big]\Big)\Big) d B_{t}+d R_{t}^{\theta} \\
&~~~~=\Big(\hat{\sigma}_{x}(t) y_{t}^{\theta}+\hat{\sigma}_{y}(t)\big(H_{t}(\theta)+\hat{R}_{t}^{\theta}(y_{t}^{\theta})\big)\Big) d B_{t}+d R_{t}^{\theta}, \\
&y_{0}^{\theta}=0 .
\end{aligned}\right.
\end{equation*}
Consequently, from (\ref{ER}) and (\ref{RR}), for $p\geq 2$,
$$\hat{\mathbb{E}}\big[\sup _{s \in[0, t]}| y_{s} ^{\theta}|^{p}\big] \leq C_p\int_{0}^{t}\Big(\hat{\mathbb{E}}\big[|y_{s}^{\theta}|^{p}\big]+|H_{s}(\theta)|^{p}+(C \theta)^{p}\Big) d s+C_{p} \theta^{p},~ t \in[0, T], ~\theta>0,$$
and, thanks to Gronwall's inequality, $\displaystyle\hat{\mathbb{E}}\Big[\sup _{t \in[0, T]}| y_{t}^ {\theta}|^{p}\Big] \leq C_p \Big(\theta^p+\int_0^T |H_{s}(\theta)|^{p}ds\Big),~\theta>0.$

Note that, for $0<\theta \leq 1$, as $F_{s}(\cdot)$ is convex, we have
\begin{equation*}
\begin{aligned}
0& \leq H_{s}(\theta) \leq F'_{s,+}(\theta)-F'_{s,+}(0) \leq F'_{s, +}(1)-F'_{s,+}(0)\\
&= \lim_{\varepsilon\downarrow0}\frac{1}{\varepsilon}\Big( \hat{\mathbb{E}}\big[\varphi_{2}(\hat{x}_s)+\varphi'_{2}(\hat{x}_s)z_{s}+\varepsilon\big(\varphi'_{2}(\hat{x}_s)z_{s}\big)\big]-\hat{\mathbb{E}}\big[\varphi_{2}(\hat{x}_s)+\varphi'_{2}(\hat{x}_s)z_{s}\big]\Big)- \hat{\mathbb{E}}_{\{\varphi_{2}(\hat{x}_{s})\}}\big[{\varphi}'_{2}(\hat{x}_s) z_{s}\big]\\
&\leq 2 \hat{\mathbb{E}}\big[|\varphi'_{2}(\hat{x}_s)z_{s}|\big]\leq C \hat{\mathbb{E}}\big[|z_{s}|\big] \leq C \hat{\mathbb{E}}\big[\sup_{s \in[0, T]}|z_{s}|\big]=: C^{*}<\infty.
\end{aligned}
\end{equation*}
Thus, since $0 \leq H_{s}(\theta) \downarrow 0$, as $\theta \downarrow 0$, $s \in[0, T]$, it follows from the bounded convergence theorem applied to $\int_0^T|H_s(\theta)|^pds$ that
\begin{equation}\label{5.11}
\hat{\mathbb{E}}\Big[\sup_{s \in[0, T]}|y_{s}^{\theta}|^{p}\Big] \rightarrow 0, \text{~as~} \theta \downarrow 0, ~p \geq 2.
\end{equation}
\begin{flushright}
$\square$
\end{flushright}
\begin{remark}\rm
From (\ref{5.10}), (\ref{RR}), (\ref{rho}) and (\ref{5.11}), we have
\begin{equation*}
\begin{aligned}
&\Big|\frac{1}{\theta}\Big(\hat{\mathbb{E}}\big[\varphi_{2}(x_{t}^{\theta})\big]-\hat{\mathbb{E}}\big[\varphi_{2}(\hat{x}_{t})\big]\Big)-\hat{\mathbb{E}}_{\{\varphi_{2}(\hat{x}_{t})\}}\big[\varphi'_{2}(\hat{x}_{t} )z_{t}\big]\Big|\\
=&\Big| H_{t}(\theta)+\hat{R}_{t} ^{\theta}(y_{t}^{\theta})\Big|
\leq H_{t}(\theta)+C\big(\theta+\hat{\mathbb{E}}\big[|y_{t}^{\theta}|\big]\big) \rightarrow 0, \text {~as~} \theta \downarrow 0, \text {~ i.e.,}
\end{aligned}
\end{equation*}
\begin{equation}\label{5.12}
\lim _{\theta \downarrow 0} \frac{1}{\theta}\Big(\hat{\mathbb{E}}\big[\varphi_{2}(x_{t}^{\theta})\big]-\hat{\mathbb{E}}\big[\varphi_{2}(\hat{x}_{t})\big]\Big)=\hat{\mathbb{E}}_{\{\varphi_{2}(\hat{x}_{t})\}}\big[\varphi'_{2}(\hat{x}_{t}) z_{t}\big],~ t \in[0, T].
\end{equation}
\end{remark}

\begin{lemma}\label{lem5.3}
The directional derivative of the cost functional $J$ is given by
\begin{equation*}
\begin{aligned}
\lim_{\theta\downarrow0}\frac{J(\hat{u}+\theta v)-J(\hat{u})}{\theta}
=&\ \hat{\mathbb{E}}_{\{\psi(\hat{u})\}}\Big[\hat{\Phi}_{x}(T)z_T +\hat{\Phi}_{y}(T)\hat{\mathbb{E}}_{\{\varphi_{4}(\hat{x}_T)\}}[\varphi_{4}'(\hat{x}_T)z_T]\\
&\ \ \ \ \ \ \ \ \ \ \ + \int_{0}^{T}\Big(\hat{l}_{x}(t)z_t +\hat{l}_y(t)\hat{\mathbb{E}}_{\{\varphi_5(\hat{x}_t)\}}[\varphi_5 '(\hat{x}_t)z_t]+\hat{l}_v(t)v_t\Big)dt\Big],
\end{aligned}
\end{equation*}
where $\displaystyle\psi(\hat{u})=\hat{\Phi}(T)+\int_{0}^{T}\hat{l }(t)dt=\Phi(\hat{x}_T,\hat{\mathbb{E}}[\varphi_4(\hat{x}_T)])+\int_{0}^{T}l (t,\hat{x}_t,\hat{\mathbb{E}}[\varphi_5(\hat{x}_t)],\hat{u}_t)dt$; for the other abbreviating notations, see \eqref{eq5.2'}.

\end{lemma}
\noindent\emph{Proof.}
For simplicity, but without restriction of the generality of the arguments, we suppose that $l=0$: $J(u)=\hat{\mathbb{E}}\big[\Phi\big(x_{T}^{u}, \hat{\mathbb{E}}[\varphi_{4}(x_{T}^{u})]\big)\big], ~u \in \mathcal{U}.$

From Lemma \ref{lem5.1}, for $ y_{t}^{\theta}=\frac{x_{t}^{\theta}-\hat{x}_{t}}{\theta}-z_t, ~t \in[0, T]$, we have
$$\hat{\mathbb{E}}\Big[\sup_{t\in[0,T]} |y_{t}^{\theta}|^{2}\Big] \rightarrow 0, \text{~as~} \theta \downarrow 0.$$
Hence, for $H(\theta):=\hat{\mathbb{E}}\big[\varphi_{4}(\hat{x}_T+\theta z_T)\big]$, $\theta \geq 0$,
$$
\big|\hat{\mathbb{E}}\big[\varphi_{4}(x_{T}^{\theta})\big]-H(\theta)\big|=\big|\hat{\mathbb{E}}\big[\varphi_{4}(\hat{x}_{T}+\theta z_{T}+\theta y_{T}^{\theta})\big]-\hat{\mathbb{E}}\big[\varphi_{4} (\hat{x}_{T}+\theta z_T)\big]\big| \leq C \theta \hat{\mathbb{E}}\big[|y_{T}^{\theta}|\big],
$$
and, thus, similar to the proof of (\ref{5.12}), we have
$$\lim _{0<\theta \downarrow 0} \frac{\hat{\mathbb{E}}\big[\varphi_{4}(x_{T}^{\theta})\big]-\hat{\mathbb{E}}\big[\varphi_{4} (\hat{x}_{T})\big]}{\theta}=H'_{+}(0),$$
where, thanks to Corollary \ref{coro4.5},
$
H'_{+}(0)=\hat{\mathbb{E}}_{\{\varphi_{4}(\hat{x}_{T})\}}\big[\varphi'_{4}(\hat{x}_{T}) z_{T}\big].
$

 Putting $\displaystyle r_{\theta}:=\frac{1}{\theta}\Big(\hat{\mathbb{E}}\big[\varphi_{4}(x_{T}^{\theta})\big]-\hat{\mathbb{E}}\big[\varphi_{4}(\hat{x}_{T})\big]\Big)-\hat{\mathbb{E}}_{\{\varphi_{4}(\hat{x}_{T})\}}\big[\varphi'_{4}(\hat{x}_{T}) z_{T}\big],~\theta>0$, and $\displaystyle \psi(\theta):= \hat{\mathbb{E}}\Big[\Phi\Big(\big(\hat{x}_{T},$ $ \hat{\mathbb{E}}[\varphi_{4}(\hat{x}_{T})]\big)+\theta\big(z_{T}, \hat{\mathbb{E}}_{\{\varphi_{4}(\hat{x}_{T})\}}\big[\varphi'_{4}(\hat{x}_{T}) z_{T}\big]\big)\Big)\Big],~ \theta \geq 0$, we have
$$\big|\hat{\mathbb{E}}\big[\Phi\big(x_{T}^{\theta}, \hat{\mathbb{E}}[\varphi_{4}(x_{T}^{\theta})]\big)\big]-\psi(\theta)\big| \leq C \theta\big(\hat{\mathbb{E}}[|y_{T}^{\theta}|^{2}]+r_{\theta}^{2}\big)^{\frac{1}{2}}, ~\theta>0,$$
with $\hat{\mathbb{E}}[|y_{T}^{\theta}|^{2}]+r_{\theta}^{2} \rightarrow 0$, as $\theta \downarrow 0$. Consequently,
\begin{equation*}
\begin{aligned}
&\lim_{0<\theta \downarrow 0} \frac{\hat{\mathbb{E}}\big[\Phi\big(x_{T}^{\theta}, \hat{\mathbb{E}}[\varphi_{4}(x_{T}^{\theta})]\big)\big]-\hat{\mathbb{E}}\big[\Phi\big(\hat{x}_{T},\hat{\mathbb{E}}[\varphi_{4}(\hat{x}_{T})]\big)\big]}{\theta}\\
=&\lim_{0<\theta \downarrow 0}\left(\frac{\hat{\mathbb{E}}\big[\Phi\big(x_{T}^{\theta}, \hat{\mathbb{E}}[\varphi_{4}(x_{T}^{\theta})]\big)\big]-\psi(\theta)}{\theta}+\frac{\psi(\theta)-\psi(0)}{\theta}                \right)
=\psi'_{+}(0),
\end{aligned}
\end{equation*}
and from Remark \ref{remark4.7},
$$\psi'_{+}(0)=\hat{\mathbb{E}}_{\{\psi(\hat{u})\}}\Big[(\partial_{x} \Phi)\big(\hat{x}_{T}, \hat{\mathbb{E}}[\varphi_{4}(\hat{x}_{T})]\big) z_{T}+(\partial_y \Phi)\big(\hat{x}_{T}, \hat{\mathbb{E}}[\varphi_{4}(\hat{x}_{T})]\big) \hat{\mathbb{E}}_{\{\varphi_{4}(\hat{x}_{T})\}}\big[\varphi'_{4}(\hat{x}_{T}) z_{T}\big]\Big].$$
The proof is complete.\hfill$\square$

\subsubsection{Duality}
In this section, we  consider the special  case where $\sigma$ and $b$
are independent of $y$, and   we still put $\beta=0$. More general cases can be studied with the same approach as that we develop here, but, of course, this is related with more involved computations. In the case we study here   (\ref{2}) becomes
\begin{equation}\left\{\begin{array}{l}
d x_{t}^{u}=\sigma(x_{t}^{u},  u_{t}) d B_{t}+b(x_{t}^{u}  , u_{t}) d t, ~ t \in[0, T],\\
x_{0}^{u}=x\in \mathbb{R}^n.
\end{array}\right.
\end{equation}
Concerning the  cost functional, we make the following assumption

\noindent(\textbf{A.3}) $\hat{\Phi}_y(T)\ge 0,\, \hat{l}_y(t)\ge 0,\, \, t\in [0,T],$ quasi-surely.

\smallskip

\noindent Of course, this assumption is, in particular, satisfied, if the partial derivates $\partial_y\Phi(.,.)$ and $\partial_y l(.,.,.,.)$ are everywhere non negative.

Recall from \eqref{eq3.2} that the cost functional is given by
\begin{equation}\label{eq5.14}
J(u):=\hat{\mathbb{E}}[\Phi(x_{T}^{u},\hat{\mathbb{E}}[\varphi_{4}(x_{T}^{u})])+\int_{0}^{T} l(t, x_{t}^{u}, \hat{\mathbb{E}}[\varphi_{5}(x_{t}^{u})], u_{t}) d t].
\end{equation}
Then from the optimality of $\hat{u}$, thanks to Lemma \ref{lem5.3}, with the notation $\psi(\hat{u})= \Phi(\hat{x}_T,\hat{\mathbb{E}}[\varphi_4(\hat{x}_T)])$ $+\int_{0}^{T} l(t, \hat{x}_{t}, \hat{\mathbb{E}}[\varphi_{5}(\hat{x}_{t})], \hat{u}_{t}) d t$ and those introduced in \eqref{eq5.2'} we have
\begin{align}\label{5.15}
0&\leq\lim_{\theta\downarrow0}\frac{J(\hat{u}+\theta v)-J(\hat{u})}{\theta}\notag \\
&=\ \hat{\mathbb{E}}_{\{\psi(\hat{u})\}}\Big[\hat{\Phi}_{x}(T)z_T +\hat{\Phi}_{y}(T)\hat{\mathbb{E}}_{\{\varphi_{4}(\hat{x}_T)\}}[\varphi_{4}'(\hat{x}_T)z_T]\notag\\
&\ \ \ \ \ \   + \int_{0}^{T}\Big(\hat{l}_{x}(t)z_t +\hat{l}_y(t)\hat{\mathbb{E}}_{\{\varphi_5(\hat{x}_t)\}}[\varphi_5 '(\hat{x}_t)z_t]+\hat{l}_v(t)v_t\Big)dt\Big]\allowdisplaybreaks\\
&=\!\!\!\sup_{P^1\in\mathcal{P}_{\{\psi(\hat{u}) \}}}\!\!\! E_{P^1}\Big[\hat{\Phi}_{x}(T)z_T +\hat{\Phi}_{y}(T)\!\!\!\sup_{P^2\in\mathcal{P}_{\{\varphi_{4}(\hat{x}_T)\}}}\!\!\!E_{P^2}[\varphi_{4}'(\hat{x}_T)z_T] \notag\\
&\ \ \ \ \ \   + \int_{0}^{T}\Big(\hat{l}_{x}(t)z_t +\hat{l}_y(t)\!\!\!\sup_{P^3\in\mathcal{P}_{\{\varphi_{5}(\hat{x}_t)\}}}\!\!\!E_{P^3}[\varphi_5 '(\hat{x}_t)z_t]+\hat{l}_v(t)v_t\Big)dt\Big]\allowdisplaybreaks\notag\\
&=\!\!\!\sup_{P^1\in\mathcal{P}_{\{\psi(\hat{u}) \}}}\!\!\! \Big\{E_{P^1}\Big[\hat{\Phi}_{x}(T)z_T + \int_{0}^{T}\Big(\hat{l}_{x}(t)z_t +\hat{l}_v(t)v_t\Big)dt\Big]\notag\\
&\ \ \ \ \ \   +E_{P^1}[\hat{\Phi}_{y}(T)]\!\!\!\sup_{P^2\in\mathcal{P}_{\{\varphi_{4}(\hat{x}_T)\}}}\hskip -0.6cm E_{P^2}[\varphi_{4}'(\hat{x}_T)z_T]+ \int_{0}^{T}\hskip -0.25cm E_{P^1}[\hat{l}_y(t)]\!\!\!\sup_{P^3\in\mathcal{P}_{\{\varphi_{5}(\hat{x}_t)\}}}\hskip -0.6cm E_{P^3}[\varphi_5 '(\hat{x}_t)z_t]dt\Big\}.\allowdisplaybreaks\notag
\end{align}
Let us now define
\begin{equation*}\label{R1}
\begin{array}{ll}
\mathcal{R}_{\{\varphi_5(\hat{x})\}}:&\hskip -0.3cm=\big\{R=(R_t):[0,T]\rightarrow\mathcal{P}\mbox{ Borel measurable} : R_t\in\mathcal{P}_{\{\varphi_5(\hat{x}_t)\}},\, t\in[0,T]\big\},\\
\mathcal{R}_{\{\varphi_5(\hat{x})|\varphi'_5(\hat{x})z\}}:&\hskip -0.3cm=\big\{R=(R_t):[0,T]\rightarrow\mathcal{P}\mbox{ Borel measurable} : R_t\in\mathcal{P}_{\{\varphi_5(\hat{x}_t)|\varphi'_5(\hat{x}_t)z_t\}},\, t\in[0,T]\big\},
\end{array}
\end{equation*}
where
\begin{equation*}
\mathcal{P}_{\{\varphi_5(\hat{x}_t)|\varphi'_5(\hat{x}_t)z_t\}}: =\big\{R\in\mathcal{P}_{\{\varphi_5(\hat{x})\}}:E_R[\varphi'_5(\hat{x}_t)z_t]=\hat{\mathbb{E}}_{\{\varphi_5(\hat{x}_t)\}}[\varphi'_5(\hat{x}_t)z_t]\big\}\subset\mathcal{P}
\end{equation*}
(cf. Definition \eqref{A.2} in Appendix 2). Here $[0,T]$ and $(\mathcal{P},d)$ are endowed with their Borel $\sigma$-algebras. Recall that $d$ is the L\'evy-Prokhorov metric on $\mathcal{P}$. From Theorem \ref{th'A.1} (a measurable selection theorem) we know that $\mathcal{R}_{\{\varphi_5(\hat{x})|\varphi'_5(\hat{x})z\}}\not=\emptyset,$ and so $\mathcal{R}_{\{\varphi_5(\hat{x})\}}\supset\mathcal{R}_{\{\varphi_5(\hat{x})|\varphi'_5(\hat{x})z\}}\not=\emptyset.$
Moreover, we observe that, for all $R=(R_t)\in\mathcal{R}_{\{\varphi_5(\hat{x})\}}$,
\begin{equation}\label{eq5.21'}
\displaystyle\int_0^TE_{P^1}[\hat{l}_y(t)]E_{R_t}[\varphi'_5(\hat{x}_t)z_t]dt\le \int_{0}^{T}\hskip -0.25cm E_{P^1}[\hat{l}_y(t)]\!\!\!\sup_{P^3\in\mathcal{P}_{\{\varphi_{5}(\hat{x}_t)\}}}\hskip -0.6cm E_{P^3}[\varphi_5 '(\hat{x}_t)z_t]dt,
\end{equation}
and, if  $R=(R_t)\in\mathcal{R}_{\{\varphi_5(\hat{x})|\varphi'_5(\hat{x})z\}}$, we have equality in \eqref{eq5.21'}.
Consequently,
\begin{equation*}\label{R1'}
\begin{array}{ll}
\displaystyle\sup_{R\in\mathcal{R}_{\{\varphi_5(\hat{x})\}} }\int_0^TE_{P^1}[\hat{l}_y(t)]E_{R_t}[\varphi'_5(\hat{x}_t)z_t]dt=\int_{0}^{T}\hskip -0.25cm E_{P^1}[\hat{l}_y(t)]\!\!\!\sup_{P^3\in\mathcal{P}_{\{\varphi_{5}(\hat{x}_t)\}}}\hskip -0.6cm E_{P^3}[\varphi_5 '(\hat{x}_t)z_t]dt,
\end{array}
\end{equation*}
and since $E_{P^1}[\hat{\Phi}_{y}(T)]\ge 0$ and $E_{P^1}[\hat{l}_y(t)]\ge 0,\, t\in[0,T],$ using the notation

\centerline{$\mathcal{P}\{\hat{u}\}:=\mathcal{P}_{\{\psi(\hat{u})\}}\times \mathcal{P}_{\{\varphi_4(\hat{x}_T)\}}\times \mathcal{R}_{\{\varphi_5(\hat{x})\}},$}

\noindent (Observe that this set does not depend on the perturbing control $u=(u_t)$)
\noindent we obtain from \eqref{5.15}
\begin{align}\label{5.15'}
0&\leq\lim_{\theta\downarrow0}\frac{J(\hat{u}+\theta v)-J(\hat{u})}{\theta}\notag \\
&=\!\!\!\sup_{(P,Q,R)\in\mathcal{P}_{\{\hat{u}\}}}\!\!\! \Big\{E_{P}\Big[\hat{\Phi}_{x}(T)z_T + \int_{0}^{T}\Big(\hat{l}_{x}(t)z_t +\hat{l}_v(t)v_t\Big)dt\Big]\\
&\ \ \ \ \ \   +E_{P}[\hat{\Phi}_{y}(T)]E_{Q}[\varphi_{4}'(\hat{x}_T)z_T]+ \int_{0}^{T}E_{P}[\hat{l}_y(t)]E_{R_t}[\varphi_5 '(\hat{x}_t)z_t]dt\Big\}.\allowdisplaybreaks\notag
\end{align}
As for the special case we consider here, (\ref{5.3}) becomes
\begin{equation}\label{eq5.3sim}
\left\{\begin{aligned}
d z_{t}=&\Big(\hat{b}_{x}(t) z_{t}+\hat{b}_{v}(t) v_{t}\Big) d t
+\Big(\hat{\sigma}_{x}(t) z_{t}+\hat{\sigma}_{v}(t) v_{t}\Big) d B_{t},~t\in[0,T], \\z_{0}=&0  .
\end{aligned}\right.
\end{equation}
Relation \eqref{5.15'} brings us to introduce the following family of adjoint BSDEs (These BSDEs are classical ones, as they are considered under a linear expectation):

\noindent 1) Under $P\in\mathcal{P}_{\{\psi(\hat{u})\}}$,
\begin{equation}\label{BSDE-P}
\left\lbrace\begin{array}{lll}
dp_s(P)&=&-\big(\hat{b}_x(s)p_s(P)+\hat{l}_x(s)\big)ds-\hat{\sigma}_x(s)q_s(P)d\langle B\rangle_s+q_s(P)dB_s+dN_s(P),\\
p_T(P)&=&\hat{\Phi}_x(T),\quad s\in[0,T],\\
& & N(P)\in \mathcal{M}_P^{2,\perp}(0,T) \mbox{ with } N_0(P)=0;
\end{array}\right.
\end{equation}
2) Under $Q\in\mathcal{P}_{\{\varphi_4(\hat{x}_T)\}}$,
 \begin{equation}\label{BSDE-Q}
\left\lbrace\begin{array}{lll}
d\tilde{p}_s(Q)&=&-\hat{b}_x(s)\tilde{p}_s(Q)ds-\hat{\sigma}_x(s)\tilde{q}_s(Q)d\langle B\rangle_s+\tilde{q}_s(Q)dB_s+d\tilde{N}_s(Q),\\
\tilde{p}_T(Q)&=&\varphi'_4(\hat{x}_T),\quad s\in[0,T],\\
& & \tilde{N}(Q)\in \mathcal{M}_Q^{2,\perp}(0,T) \mbox{ with } \tilde{N}_0(Q)=0;
\end{array}\right.
\end{equation}
3) Under $R_t,\ t\in [0,T)$, for $R=(R_t)\in\mathcal{R}_{\{\varphi_5(\hat{x})\}}$,
 \begin{equation}\label{BSDE-Rt}
\left\lbrace\begin{array}{lll}
p_s(t,R_t)&=&-\hat{b}_x(s)p_s(t,R_t)ds-\hat{\sigma}_x(s)q_s(t,R_t)d\langle B\rangle_s+q_s(t,R_t)dB_s+dN_s(t,R_t),\\
p_t(t,R_t)&=&\varphi'_5(\hat{x}_t),\quad s\in[0,t],\\
& & N(t,R_t)\in \mathcal{M}_R^{2,\perp}(0,T) \mbox{ with } N_0(t,R_t)=0.
\end{array}\right.
\end{equation}
\begin{remark}\rm
1) For the above BSDEs we consider the measurable space $(\Omega,\mathcal{B}(\Omega))$ endowed with the filtration $\mathbb{F}^B=(\mathcal{F}_s)$ generated by the $G$-Brownian motion $B$ (Recall that $B$ has been introduced as coordinate process on $\Omega$). For a given probability measure $P$ over $(\Omega,\mathcal{B}(\Omega))$ the associated
filtration is the one augmented by all $P$-null sets: $\mathbb{F}^P=\mathbb{F}^B\vee \mathcal{N}_P$.

2) Note that, under any $P\in\mathcal{P}$, the $G$-Brownian motion $B$ is only a continuous square integrable martingale, and so the martingale representation  may not hold for $(B,\mathbb{F}^P)$. So it is necessary to introduce the second square integrable $P$-martingale $N(P)$ with $N_0(P)=0$ and joint quadratic variation $\langle B, N(P)\rangle^{P}\big(=(\langle B, N(P)\rangle^{P}_s)\big)=0$ (We write $N(P)\in \mathcal{M}_P^{2,\perp}(0,T)$).

3) Recall that $\langle B\rangle$ is the quadratic variation process of the $G$-Brownian motion $B$ under $\hat{\mathbb{E}}$: For all $\pi_{t}^{N}=\{0=t_0^N<t_1^N<\dots<t_N^N=t\}, ~N \geq 1$, sequence of partitions of $[0, t]$ with mesh $|\pi_{t}^{N}|=\displaystyle\max_{0\le j\le N-1}(t_{j+1}^N-t_j^N) \rightarrow 0~(N \rightarrow \infty)$,
$$\hat{\mathbb{E}}\bigg[\bigg|\sum_{j=0}^{N-1}(B_{t_{j+1}^{N}}-B_{t_{j}^{N}})^{2}-\langle B\rangle_{t}\bigg|^{2}\bigg] \rightarrow 0~(N \rightarrow \infty).$$
And so, for all $P \in \mathcal{P}$, $\langle B\rangle$ coincides $P$-a.s. with the quadratic variation process $\langle B\rangle^{P}$ of $B$ as $P$-martingale, $\langle B\rangle_{t}^{P}=\langle B\rangle_{t}, ~t \in[0, T], ~P\text{-a.s.}$ Also recall that, under the $G$-expectation the increments of $\langle B\rangle$ are independent and stationary, and $\underline{\sigma}^2ds\le d\langle B\rangle_s\le \overline{\sigma}^2ds,$ $ds$-a.e., quasi-surely.

\end{remark}

Following El Karoui and Huang \cite{22} and Buckdahn et al. \cite{23}, we see that, for all $P\in\mathcal{P}$, there exists a unique triplet of processes $(p(P), q(P), N(P)) \in M_{P}^{2}(0, T) \times M_{P}^{2}(0, T) \times \mathcal{M}_{P}^{2, \perp}(0, T )$ which solves the adjoint equations (\ref{BSDE-P}) and (\ref{BSDE-Q}) (equation (\ref{BSDE-Q}) with $Q$ instead of $P$), respectively. The same we also have for the BSDE \eqref{BSDE-Rt}, only that here the BSDE is considered over the time interval $[0,t]$, so that the unique solution triplet $(p(t,R_t), q(t,R_t), N(t,R_t))$ belongs to $M_{R}^{2}(0, t) \times M_{R}^{2}(0, t) \times \mathcal{M}_{R}^{2, \perp}(0, t)$, $t\in[0,T)$. Moreover, standard BSDE estimates using that the coefficients $\hat{b}_x,\, \hat{\sigma}_x,\ \hat{l}_x$ are bounded, show that, for all $p\ge 1$, there is some constant $C_p\in\mathbb{R}_+$ (independent of the underlying probability measure $P\in\mathcal{P}$) s.t.
\begin{equation}\label{B1}
\displaystyle E_P\big[\sup_{s\in[0,T]}|p_s(P)|^p+\big(\int_0^T|q_s(P)|^2d\langle B\rangle_s+\langle N(P)\rangle_T \big)^{p/2}\big]\le C_p.
\end{equation}
Similar estimates we have for the solution $(p(t,R_t), q(t,R_t), N(t,R_t))\in M_{R}^{2}(0, t) \times M_{R}^{2}(0, t) \times \mathcal{M}_{R}^{2, \perp}(0, t)$ of BSDE \eqref{BSDE-Rt}, for all $t\in[0,T]$, only that unlike in \eqref{B1}, here $T$ has to be replaced by $t$. The constant $C_p$ in the estimate of $(p(t,R_t), q(t,R_t), N(t,R_t))$ is again independent of $R=(R_t)\in\mathcal{R}_{\{\varphi_5(\hat{x})\}}$ but also independent of $t\in[0,T].$

Applying now It\^{o}'s formula to $p_s(P) z_s$, we have
\begin{equation}
\begin{aligned}
d(p_s(P)z_s)=\big(p_s(P)\hat{b}_{v}(s) v_{s}-\hat{l}_x(s)z_s\big) d s+\zeta_s(P)d B_{s}+q_{s}(P)\hat{\sigma}_{v}(s)v_s  d\langle B\rangle_{s}+z_{s} d N_{s}(P),
\end{aligned}
\end{equation}
where $\zeta_s(P):=p_{s}(P)(\hat{\sigma}_{x}(s) z_{s}+\hat{\sigma}_{v}(s)v_{s}) +z_{s} q_{s}(P)$.
As $z_0=0$,
\begin{equation}
\begin{array}{lll}
p_T(P)z_T&=&\displaystyle\int_0^T \big(p_{s}(P)\hat{b}_{v}(s) v_{s}-\hat{l}_x(s)z_s\big) d s\\
& &\displaystyle+\int_0^T\zeta_s(P)d B_{s}+\int_0^Tq_{s}(P)\hat{\sigma}_{v}(s)v_s d\langle B\rangle_{s}+\int_0^Tz_{s} d N_{s}(P),
\end{array}
\end{equation}
where $\displaystyle\int_0^{\cdot}\zeta_s(P)d B_{s}$ and $\displaystyle\int_0^{\cdot}z_{s} d N_{s}(P)$ are $P$-martingales. Indeed, from our estimates it follows that
$$
\begin{aligned}
E_{P}\Big[\Big(\int_{0}^{T}|z_{t}|^{2} d\langle N(P)\rangle_t\Big)^{\frac{1}{2}}\Big]& \leq E_{P}\big[\sup_{0 \leq t \leq T} |z_{t}|\langle N(P) \rangle_{T}^{\frac{1}{2}}\big] \\
&\leq \Big(\hat{\mathbb{E}}\big[\sup _{0 \leq t \leq T}|z_{t}|^{2}\big]\Big)^{\frac{1}{2}}\Big(E_{P}\big[\langle N(P)\rangle_{T}\big]\Big)^{\frac{1}{2}}<+\infty,
\end{aligned}
$$
and with similar arguements we also see that $\displaystyle E_{P}\Big[\Big(\int_{0}^{T}|\zeta_{t}(P)|^{2} d\langle B\rangle_t\Big)^{\frac{1}{2}}\Big]<+\infty$.

Thus, recallig that $p_T(P)=\hat{\Phi}_x(T)$, we have
\begin{equation}\label{5.22}
\begin{array}{ll}
&\displaystyle E_{P}\big[\hat{\Phi}_{x}(T)z_T+\int_0^T\big(\hat{l}_x(s)z_s+\hat{l}_v(s)v_s\big)ds\big]\\
&=\displaystyle E_{P}\Big[\int_0^T  v_{s}\Big(\big(p_{s}(P) \hat{b}_{v}(s)+\hat{l}_v(s)\big)ds +q_{s}(P)\hat{\sigma}_{v}(s)d\langle B\rangle_s\Big)\Big].
\end{array}
\end{equation}
An analogous argument but with using now the solution of BSDE \eqref{BSDE-Q} yields, for $Q\in\mathcal{P}$,
\begin{equation}\label{5.21}
E_{Q}[\varphi_{4}'(\hat{x}_T)z_T]=E_{Q}[\tilde{p}_T(Q) z_T]=E_{Q}\Big[\int_0^T  v_{s}\Big(\tilde{p}_{s}(Q)\hat{b}_{v}(s)ds +\tilde{q}_{s}(Q) \hat{\sigma}_{v}(s)d\langle B\rangle_s\Big)\Big].
\end{equation}
Finally, making use in the same way of the solution $(p(t,R_t), q(t,R_t), N(t,R_t))\in M_{R}^{2}(0, t) \times M_{R}^{2}(0, t) \times \mathcal{M}_{R}^{2, \perp}(0, t)$ of BSDE \eqref{BSDE-Rt}, we obtain, for $t\in[0,T],$
\begin{equation}\label{5.21-Rt}
E_{R_t}[\varphi_{5}'(\hat{x}_t)z_t]=E_{R_t}[p_t(t,R_t) z_t]=E_{R_t}\Big[\int_0^t  v_{s}\Big(p_{s}(t,R_t)\hat{b}_{v}(s)ds +q_{s}(t,R_t) \hat{\sigma}_{v}(s)d\langle B\rangle_s\Big)\Big].
\end{equation}

Let us introduce now
\begin{equation}\label{Theta1}
\begin{array}{lll}
 \Theta[P,Q,R](v)&=&\displaystyle E_{P}\big[\hat{\Phi}_{x}(T)z_T+\int_0^T\big(\hat{l}_x(s)z_s+\hat{l}_v(s)v_s\big)ds\big]\\
& &+\displaystyle E_{P}\big[\hat{\Phi}_y(T)\big]E_{Q}\big[\varphi_{4}'(\hat{x}_T)z_T\big]+\displaystyle\int_0^T  E_{P}\big[\hat{l}_y(t)\big] E_{R_t}\big[\varphi_{5}'(\hat{x}_t)z_t\big],
\end{array}
\end{equation}
and from the above computation we see that
\begin{equation}\label{Theta1}
\begin{array}{lll}
\hskip -0.45cm\Theta[P,Q,R](v)\hskip -0.25cm&=&\hskip -0.25cm\displaystyle E_{P}\Big[\int_0^T  v_{s}\Big(\big(p_{s}(P) \hat{b}_{v}(s)+\hat{l}_v(s)\big)ds +q_{s}(P)\hat{\sigma}_{v}(s)d\langle B\rangle_s\Big)\Big]\\
& \hskip -0.25cm&+\displaystyle E_{P}\Big[\hat{\Phi}_y(T)\Big]E_{Q}\Big[\int_0^T  v_{s}\Big(\tilde{p}_{s}(Q)\hat{b}_{v}(s)ds +\tilde{q}_{s}(Q) \hat{\sigma}_{v}(s)d\langle B\rangle_s\Big)\Big]\\
& &\hskip -0.25cm+\displaystyle\int_0^T  E_{P}\Big[\hat{l}_y(t)\Big] E_{R_t}\Big[\int_0^t  v_{s}\Big(p_{s}(t,R_t)\hat{b}_{v}(s)ds +q_{s}(t,R_t) \hat{\sigma}_{v}(s)d\langle B\rangle_s\Big)\Big]dt.
\end{array}
\end{equation}
In order to give to \eqref{Theta1}  another form, we make the convention that $p_{s}(t,R_t):=0, q_{s}(t,R_t): =0,$ for $t<s\le T$, and we define the probability measure $\widetilde{R}:=\displaystyle\int_0^T\frac{1}{T}dt\cdot\big(\delta_t\otimes R_t\big)$ over the probability space $([0,T]\times\Omega,\mathcal{B}([0,T])\otimes\mathcal{F})$. Here $\delta_t$ denotes the Dirac measure over $[0,T]$ with mass at $t$. Then, with $(t,\omega)\mapsto \big(p_{s}(t,R_t)(\omega), q_{s}(t,R_t)(\omega)\big)$ and $t\mapsto E_{P}\Big[\hat{l}_y(t)\Big]$ interpreted as random variables over $[0,T]\times\Omega$, we have
\begin{equation}\label{Theta2}
\begin{array}{lll}
& &\displaystyle\int_0^T  E_{P}\Big[\hat{l}_y(t)\Big] E_{R_t}\Big[\int_0^t  v_{s}\Big(p_{s}(t,R_t)\hat{b}_{v}(s)ds +q_{s}(t,R_t) \hat{\sigma}_{v}(s)d\langle B\rangle_s\Big)\Big]dt\\
&=&\displaystyle TE_{\widetilde{R}}\Big[\int_0^Tv_s\Big(E_P[\hat{l}_y(\cdot)]\Big(p_{s}(\cdot,R_.)\hat{b}_{v}(s)ds +q_{s}(\cdot,R_.) \hat{\sigma}_{v}(s)d\langle B\rangle_s\Big)\Big)\Big].
\end{array}
\end{equation}

\vskip 1cm

Let us define $\Omega_{\{T\}}:=[0,T]\times\Omega$ and embed the probabilites $P$ and $Q$ in the space of probabilities over $(\Omega_{\{T\}},\mathcal{B}([0,T])\otimes\mathcal{F})$ in a canonical way by making the identification $P:=\delta_T\otimes P$ and $Q:=\delta_T\otimes Q$. Then, thanks to \eqref{Theta1} and \eqref{Theta2},
\begin{equation}\label{Theta3}
\begin{array}{ll}
&\Theta[P,Q,R](v)\\
&=\quad\displaystyle \int_{\Omega_{\{T\}}}\int_0^Tv_s \Big\{\Big(\big(p_s(P)\hat{b}_v(s)+\hat{l}_v(s)\big)ds+q_s(P)\hat{\sigma}_v(s)d\langle B\rangle_s\Big)dP\\
&\hskip 3cm + E_P[\hat{\Phi}_y(T)]\Big(\tilde{p}_s(Q)\hat{b}_v(s)ds+\tilde{q}_s(Q)\hat{\sigma}_v(s)d\langle B\rangle_s \Big)dQ\\
&\hskip 3cm +TE_P[\hat{l}_y(\cdot)]\Big(p_s(\cdot,R_.)\hat{b}_v(s)ds+q_s(\cdot,R_.)\hat{\sigma}_v(s)d\langle B\rangle_s\Big)d\tilde{R}\Big\}\\
&=\quad\displaystyle \int_{\Omega_{\{T\}}}\int_0^Tv_s\big\{\hat{b}_v(s)dsdp_s(P,Q,R)+\hat{l}_v(s)dsdP+\hat{\sigma}_v(s)d\langle B\rangle_sdq_s(P,Q,R)\big\},
\end{array}
\end{equation}
where

\centerline{$\begin{array}{lll}
dp_s(P,Q,R):&=&p_s(P)dP+E_P[\hat{\Phi}_y(T)]\tilde{p}_s(Q)dQ+TE_p[\hat{l}_y(\cdot)]p_s(\cdot,R_.) d\tilde{R},\\
dq_s(P,Q,R):&=&q_s(P)dP+E_P[\hat{\Phi}_y(T)]\tilde{q}_s(Q)dQ+TE_P[\hat{l}_y(\cdot)]q_s(\cdot,R_.)d\tilde{R}.
\end{array}$}

\noindent We remark that $dsdp_s(P,Q,R)$ and $d\langle B\rangle_sdq_s(P,Q,R)$ are signed measures on $\Omega_{\{T\}}\times[0,T]$ not depending on $v$ and so neither on the perturbing control $u$.
Then, from  (\ref{5.15'}), \eqref{Theta1} and \eqref{Theta3}, and with the Hamiltonian measure
$$dH_v(s,P,Q,R):=\hat{b}_v(s)dsdp_s(P,Q,R)+\hat{l}_v(s)dsdP+\hat{\sigma}_v(s)d\langle B\rangle_sdq_s(P,Q,R)$$
we have, for all $u\in\mathcal{U}$ (Recalling that $v=u-\hat{u}$) that
\begin{equation}\label{5.23}
\begin{array}{lll}
0&\le&\hskip -0.2cm\displaystyle \sup_{(P,Q,R)\in\mathcal{P}\{\hat{u}\}}\Theta[P,Q,R](u-\hat{u})\\
&=&\displaystyle \sup_{(P,Q,R)\in\mathcal{P}\{\hat{u}\}}\displaystyle \int_{\Omega_{\{T\}}}\int_0^T\hskip -0.2cm(u_s-\hat{u}_s) dH_v(s,P,Q,R).
\end{array}
\end{equation}
Observe that (\ref{5.23}) gives a necessary condition for the optimality of the control $\hat{u}\in\mathcal{U}$. We resume our main result:
\begin{theorem}\label{th5.5}
 Suppose $(\textbf{\rm A.1})$-$(\textbf{\rm A.3})$ where $b$, $\sigma$ are independent of $y$, and let $\hat{u}$ be an optimal control with state trajectory $\hat{x}=(\hat{x}_{t})$. Then (\ref{5.23}) gives a necessary optimality condition satisfied by  all $u\in\mathcal{U}$.
\end{theorem}

In the particular case when $\hat{l}_y(t)=0,$ quasi-surely, $dt$-a.s., and $\hat{\Phi}_y(T)$ is deterministic, by using an argument developed by Hu and Ji \cite{14} based on Sion's minimax theorem, we can simplify  the necessary optimality condition \eqref{5.23}.
Indeed, let us suppose

\smallskip

\noindent(\textbf{A.3'}) $l(t,x,y,u)=l(t,x,u),\, \Phi(x,y)=\Phi_1(x)+\Phi_2(y),\, (t,x,y,u)\in[0,T]\times \mathbb{R}\times\mathbb{R}\times U.$

\smallskip

\noindent We observe that under Assumption (A.3') $\hat{l}_y(t)=0$ everywhere on $[0,T]\times \Omega$ and $\hat{\Phi}_y(T)=(\Phi_2)_y(\hat{\mathbb{E}}[\varphi_4(t)])$ is deterministic. Then \eqref{5.15'} takes the simpler form
\begin{equation}\label{5.15''}
\begin{array}{ll}
0&\leq\displaystyle\lim_{\theta\downarrow0}\frac{J(\hat{u}+\theta v)-J(\hat{u})}{\theta}\\
&=\, \, \!\displaystyle\sup_{(P,Q)\in\mathcal{P}_{\{\psi(\hat{u})\}}\times \mathcal{P}_{\{\varphi_4(\hat{x}_T)\}}}\!\!\! \Big\{E_{P}\Big[\hat{\Phi}_{x}(T)z_T + \int_{0}^{T}\hskip -0.3cm \Big(\hat{l}_{x}(t)z_t +\hat{l}_v(t)v_t\Big)dt\Big]+\hat{\Phi}_{y}(T)E_{Q}[\varphi_{4}'(\hat{x}_T)z_T]\Big\}.
\end{array}
\end{equation}
We remark that the function $F:\big(\mathcal{P}_{\{\psi(\hat{u})\}}\times \mathcal{P}_{\{\varphi_4(\hat{x}_T)\}}\big)\times \mathcal{U}\rightarrow \mathbb{R}$, defined by
\begin{align}\label{5.15a'}
&F\big((P,Q),u \big)\notag\\
=&E_{P}\Big[\hat{\Phi}_{x}(T)z_T^{u} + \int_{0}^{T}\Big(\hat{l}_{x}(t)z_t^{u} +\hat{l}_v(t)(u_t-\hat{u}_t\Big)dt\Big]+\displaystyle \hat{\Phi}_{y}(T)E_{Q}[\varphi_{4}'(\hat{x}_T)z_T^{u}],\notag\\
&\hskip 5cm\big((P,Q),u\big)\in \big(\mathcal{P}_{\{\psi(\hat{u})\}}\times \mathcal{P}_{\{\varphi_4(\hat{x}_T)\}}\big)\times \mathcal{U},\allowdisplaybreaks\notag
\end{align}
is affine in $(P,Q)$ over $\mathcal{P}_{\{\psi(\hat{u})\}}\times \mathcal{P}_{\{\varphi_4(\hat{x}_T)\}}$ and affine in $u$ over $\mathcal{U}$ (Recall SDE \eqref{eq5.3sim} for $z^{u}=z$):
\begin{align}
&F(\lambda (P,Q)+(1-\lambda)(P',Q'),u)=\lambda F((P,Q),u)+(1-\lambda)F((P',Q'),u),\notag\\
&F((P,Q),\lambda u+(1-\lambda)u')=\lambda F((P,Q), u)+(1-\lambda)F((P,Q),u'),\notag\\
& \hskip 2cm  (P',Q'),(P,Q)\in \mathcal{P}_{\{\psi(\hat{u})\}}\times \mathcal{P}_{\{\varphi_4(\hat{x}_T)\}},\, u,u'\in \mathcal{U},\, \lambda\in[0,1].
\end{align}
The fact that $\mathcal{P}_{\{\psi(\hat{u})\}}\times \mathcal{P}_{\{\varphi_4(\hat{x}_T)\}}$ is a non-void convex and weakly compact subset of a linear topological space (that of the pairs of bounded signed measures) and $\mathcal{U}$ is a convex subset (Recall that the control state space $U$ is convex) of a linear topological space, Sion's minimax theorem applies,
\begin{equation*}
\begin{array}{lll}
0&\le&\displaystyle\inf_{u\in\mathcal{U}}\,\,\!\!\sup_{(P,Q)} \Big\{E_{P}\Big[\hat{\Phi}_{x}(T)z_T + \int_{0}^{T}\Big(\hat{l}_{x}(t)z_t +\hat{l}_v(t)(u_t-\hat{u}_t)\Big)dt\Big]+\displaystyle \hat{\Phi}_{y}(T)E_{Q}[\varphi_{4}'(\hat{x}_T)z_T]\Big\}\\
&=&\displaystyle\sup_{(P,Q)} \inf_{u\in\mathcal{U}}\Big\{E_{P}\Big[\hat{\Phi}_{x}(T)z_T + \int_{0}^{T}\Big(\hat{l}_{x}(t)z_t +\hat{l}_v(t)(u_t-\hat{u}_t)\Big)dt\Big]
+\displaystyle \hat{\Phi}_{y}(T) E_{Q}[\varphi_{4}'(\hat{x}_T)z_T]\Big\},
\end{array}
\end{equation*}
where the supremum is taken over all $(P,Q)\in\mathcal{P}_{\{\psi(\hat{u})\}}\times \mathcal{P}_{\{\varphi_4(\hat{x}_T)\}}.$
By using the weak compactness of $\mathcal{P}_{\{\psi(\hat{u})\}}\times \mathcal{P}_{\{\varphi_4(\hat{x}_T)\}}$ a standard argument allows to show that there exists $(P^*,Q^*)\in\mathcal{P}_{\{\psi(\hat{u})\}}\times \mathcal{P}_{\{\varphi_4(\hat{x}_T)\}}$ for which the latter supremum is attained (see also \cite{14}, proof of Theorem 4.6), i.e.,
\begin{equation}\label{aaa}
\begin{array}{lll}
0&\le&\displaystyle \inf_{u\in\mathcal{U}}\Big\{E_{P^*}\Big[\hat{\Phi}_{x}(T)z_T + \int_{0}^{T}\Big(\hat{l}_{x}(t)z_t +\hat{l}_v(t)(u_t-\hat{u}_t)\Big)dt\Big]\\
& &\hskip 6.5cm +\displaystyle \hat{\Phi}_{y}(T) E_{Q^*}[\varphi_{4}'(\hat{x}_T)z_T]\Big\}.
\end{array}
\end{equation}
This makes that we only have to use the adjoint BSDEs  \eqref{BSDE-P} and \eqref{BSDE-Q} under $P^*$ and $Q^*$, respectively, and the necessary optimality condition \eqref{5.23} takes the form
\begin{equation}\label{5.23c}
\begin{array}{rll}
0&\le &\displaystyle\int_\Omega\int_0^T\hskip -0.3cm (u_s-\hat{u}_s)\Big\{\hat{b}_v(s)dsdp_s(P^*,Q^*)+\hat{l}_v(s)dsdP^*+\hat{\sigma}_v(s)d\langle B\rangle_sdq_s(P^*,Q^*)\Big\},\, u\in\mathcal{U},
\end{array}
\end{equation}
where
\begin{equation*}\begin{array}{lll}
dp_s(P^*,Q^*):&=&p_s(P^*)dP^*+\hat{\Phi}_y(T)\tilde{p}_s(Q^*)dQ^*,\\
dq_s(P^*,Q^*):&=&q_s(P^*)dP^*+\hat{\Phi}_y(T)\tilde{q}_s(Q^*)dQ^*.
\end{array}
\end{equation*}

Finally, from the arbitrariness of $u\in\mathcal{U}$ we obtain
\begin{theorem}\label{th5.6}
Suppose $(\textbf{\rm A.1})$, $(\textbf{\rm A.2})$ and $(\textbf{\rm A.3'})$ where $b$ are $\sigma$ do not depend on $y$, and let $\hat{u}$ be an optimal control with the associated state trajectory $\hat{x}=(\hat{x}_{t})$. Then there exists $(P^*,Q^*)\in \mathcal{P}_{\{\psi(\hat{u})\}}\times \mathcal{P}_{\{\varphi_4(\hat{x}_T)\}}$ such that,  all $u\in\mathcal{U},$
\begin{equation}\label{5.23d}
\begin{array}{lll}
0&\le &\displaystyle (u_s-\hat{u}_s)\Big\{ \hat{b}_v(s)dsdp_s(P^*,Q^*)+\hat{l}_v(s)dsdP^*+\hat{\sigma}_v(s)d\langle B\rangle_sdq_s(P^*,Q^*)\Big\}.
\end{array}
\end{equation}

\end{theorem}

\subsection{ Sufficient conditions for optimality}
In this section, we continue to consider the case discussed in Section 5.1.2. We define the Hamiltonian random field
$$dH(t,x,u,p,q):=H_1(x,u,p)dt+H_2(x,u,q)d\langle B\rangle_t,$$
with $H_1(x,u,p):=b(x,u) p$ and $H_2(x,u,q):=\sigma(x, u) q,$ we make the following additional assumption:

\noindent(\textbf{A.4}) The function $\Phi$ is convex in $(x,y)$; the running cost $l(t,.,.,.)$ is convex, for all $t\in[0,T]$; the functions $\varphi_4$ and $\varphi_5$ are convex;  the Hamiltonian random field $dH(t,x,u,p,q)$ is  convex in $(x, u)$ (defined by the convexity of $H_1(\cdot,\cdot,p)$ and that of $H_2(\cdot,\cdot,q)$).

\begin{theorem}
Assume the conditions $(\textbf{\rm A.1})$-$(\textbf{\rm A.4})$ are satisfied and let $\hat{u} \in \mathcal{U}$ be a control process with associated state process  $\hat{x}=(\hat{x}_{t})$, and let $(p(P), q(P),N(P))$, $(\tilde{p}(Q), \tilde{q}(Q),\tilde{N}(Q))$ and $(p(t,R_t),q(t,R_t),N(t,R_t)),\, t\in[0,T],\, (P,Q,R)\in\mathcal{P}\{\hat{u}\}$, be the solution of BSDE
\eqref{BSDE-P}, \eqref{BSDE-Q} and \eqref{BSDE-Rt}, respectively. If (\ref{5.23}) holds for  all $u\in\mathcal{U}$, then $\hat{u}$ is an optimal control.

\end{theorem}
\noindent\emph{Proof.}
Let $u\in\mathcal{U}$ be any admissible control. From (\ref{eq5.14}), with
\begin{equation*}\begin{array}{lll}
\xi^{u}&:=&\displaystyle \Phi(x_{T}^{u},\hat{\mathbb{E}}[\varphi_{4}(x_{T}^{u})])+\int_{0}^{T} l(t, x_{t}^{u}, \hat{\mathbb{E}}[\varphi_{5}(x_{t}^{u})], u_{t}) d t,\\
\hat{\xi}&:=&\displaystyle \Phi(\hat{x}_{T},\hat{\mathbb{E}}[\varphi_{4}(\hat{x}_{T})])+\int_{0}^{T} l(t, \hat{x}_{t}, \hat{\mathbb{E}}[\varphi_{5}(\hat{x}_{t})], \hat{u}_{t})dt,
\end{array}
\end{equation*}
we have $J(u)-J(\hat{u})=\hat{\mathbb{E}}[\xi]-\hat{\mathbb{E}}[\hat{\xi}].$

Since the function $F(\lambda):=\hat{\mathbb{E}}[\hat{\xi}+\lambda(\xi^{u}-\hat{\xi})], ~\lambda \in[0,1],$ is convex,
$F(1)-F(0) \geq F'_{+}(0)$.  Thus, from Lemma \ref{lem4.4} we have
\begin{equation}\label{5.24}
\hat{\mathbb{E}}[\xi^{u}]-\hat{\mathbb{E}}[\hat{\xi}] \geq\lim_{\lambda\searrow 0} \frac{1}{\lambda}\big(\hat{\mathbb{E}}\big[\hat{\xi}+\lambda(\xi^{u}-\hat{\xi}\,)\big]-\hat{\mathbb{E}}\big[\hat{\xi}\,\big]\big)=\hat{\mathbb{E}}_{\{\hat{\xi}\}}[\xi^{u}-\hat{\xi}\,].
\end{equation}
On the other hand, from the convexity of $\Phi$, we get
$$\begin{aligned}
&\Phi\big(x_{T}^{u}, \hat{\mathbb{E}}[\varphi_{4}(x_{T}^{u})]\big)-\Phi\big(\hat{x}_{T}, \hat{\mathbb{E}}[\varphi_{4}(\hat{x}_{T})]\big)\\
&\geq \Phi_{x}\big(\hat{x}_{T}, \hat{\mathbb{E}}[\varphi_{4}(\hat{x}_{T})]\big)(x_{T}^{u}-\hat{x}_{T})+\Phi_{y}\big(\hat{x}_{T}, \hat{\mathbb{E}}[\varphi_{4}(\hat{x}_{T})]\big)\big(\hat{\mathbb{E}}[\varphi_{4}(x_{T}^{u})]-\hat{\mathbb{E}}[\varphi_{4}(\hat{x}_{T})]\big).
\end{aligned}
$$
Using now the convexity of $\varphi_4$ as well as (\ref{5.24}), but now with $\xi^{u}=\varphi_{4}(x_{T}^{u})$and $\hat{\xi}=\varphi_{4}(\hat{x}_{T})$, we see that
$$\begin{aligned}
&\hat{\mathbb{E}}[\varphi_{4}(x_{T}^{u})]-\hat{\mathbb{E}}[\varphi_{4}(\hat{x}_{T})]\ge \hat{\mathbb{E}}_{\{\varphi_{4}(\hat{x}_{T})\}}[\varphi_{4}(x_{T}^{u})-\varphi_{4}(\hat{x}_{T})]\\
&\ge \hat{\mathbb{E}}_{\{\varphi_{4}(\hat{x}_{T})\}}\big[\varphi'_{4}(\hat{x}_{T})\big(x_{T}^{u}-\hat{x}_{T}\big)\big],
\end{aligned}
$$
and from the non negativity of $\hat{\Phi}_y(T)$ we obtain
$$\begin{aligned}
&\Phi\big(x_{T}^{u}, \hat{\mathbb{E}}[\varphi_{4}(x_{T}^{u})]\big)-\Phi\big(\hat{x}_{T}, \hat{\mathbb{E}}[\varphi_{4}(\hat{x}_{T})]\big)\\
&\geq \Phi_{x}\big(\hat{x}_{T}, \hat{\mathbb{E}}[\varphi_{4}(\hat{x}_{T})]\big)(x_{T}^{u}-\hat{x}_{T})+\Phi_{y}\big(\hat{x}_{T}, \hat{\mathbb{E}}[\varphi_{4}(\hat{x}_{T})]\big)\hat{\mathbb{E}}_{\{\varphi_{4}(\hat{x}_{T})\}}\big[\varphi'_{4}(\hat{x}_{T})\big(x_{T}^{u}-\hat{x}_{T}\big)\big]\\
&= \hat{\Phi}_{x}(T)(x_{T}^{u}-\hat{x}_{T})+\hat{\Phi}_{y}(T)\hat{\mathbb{E}}_{\{\varphi_{4}(\hat{x}_{T})\}}\big[\varphi'_{4}(\hat{x}_{T})\big(x_{T}^{u}-\hat{x}_{T}\big)\big]
\end{aligned}
$$
(Recall the notations introduced in \eqref{eq5.2'}). Similarly, we see that, thanks to the convexity of $l(t,.,.,.)$ and $\varphi_5$ as well as the non negativity of $\hat{l}_y(T)$,
$$\begin{aligned}
& l(t, x_{t}^{u}, \hat{\mathbb{E}}[\varphi_{5}(x_{t}^{u})], u_{t})-l(t, \hat{x}_{t}, \hat{\mathbb{E}}[\varphi_{5}(\hat{x}_{t})], \hat{u}_{t})\\
&\geq \hat{l}_{x}(t)(x_{t}^{u}-\hat{x}_{t})+\hat{l}_{y}(t)\hat{\mathbb{E}}_{\{\varphi_{5}(\hat{x}_{t})\}}\big[\varphi'_{5}(\hat{x}_{t})\big(x_{t}^{u}-\hat{x}_{t}\big)\big]+\hat{l}_v(t)(u_t-\hat{u}_t).
\end{aligned}
$$
Hence, with the notation $\mathcal{P}_{\{\hat{\xi}\}}=\mathcal{P}_{\{\psi(\hat{u})\}}$ ($\psi(\hat{u})$ has been introduced in Lemma \ref{lem5.3}), by summarising the above computations we obtain
\begin{equation}\label{equ.5.2a}
\begin{array}{lll}
\hskip -0.5cm J(u)-J(\hat{u})\hskip -0.25cm &=&\hskip -0.25cm \hat{\mathbb{E}}[\xi]-\hat{\mathbb{E}}[\hat{\xi}]\ge\hat{\mathbb{E}}_{\{\hat{\xi}\}}[\xi^{u}-\hat{\xi}\,]\\
&\ge &\hskip -0.25cm \displaystyle\hat{\mathbb{E}}_{\{\psi(\hat{u})\}} \big[ \hat{\Phi}_{x}(T)(x_{T}^{u}-\hat{x}_{T})+\hat{\Phi}_{y}(T)\hat{\mathbb{E}}_{\{\varphi_{4}(\hat{x}_{T})\}}\Big[\varphi'_{4}(\hat{x}_{T})\big(x_{T}^{u}-\hat{x}_{T}\big)\big]\\
& & \displaystyle \hskip -0.25cm+\int_0^T\hskip -0.3cm\big(\hat{l}_{x}(t)(x_{t}^{u}-\hat{x}_{t})+\hat{l}_{y}(t)\hat{\mathbb{E}}_{\{\varphi_{5}(\hat{x}_{t})\}}\big[\varphi'_{5}(\hat{x}_{t})\big(x_{t}^{u}-\hat{x}_{t}\big)\big]+\hat{l}_v(t)(u_t-\hat{u}_t)
\big)dt\Big].\\
\end{array}
\end{equation}
Let us introduce now the following notations related with our Hamiltonian:
\begin{equation}\label{equ.5.2Ham}
\begin{array}{lll}
dH^{u,P}(s)&:=&b^{u}(s)p_s(P)ds+\sigma^{u}(s)q_s(P)d\langle B\rangle_s,\\
d\hat{H}^{P}(s)&:=&\hat{b}(s)p_s(P)ds+\hat{\sigma}(s)q_s(P)d\langle B\rangle_s,\\
d\hat{H}_x^{P}(s)&:=&\hat{b}_x(s)p_s(P)ds+\hat{\sigma}_x(s)q_s(P)d\langle B\rangle_s,\\
d\hat{H}_v^{P}(s)&:=&\hat{b}_v(s)p_s(P)ds+\hat{\sigma}_v(s)q_s(P)d\langle B\rangle_s,
\end{array}
\end{equation}
where $(b^{u},\sigma^{u})(s):=(b,\sigma)(x^{u}_s,u_s)$ and $l^{u}(s):=l(x^{u}_s,\hat{\mathbb{E}}[\varphi_5(x_s^{u})],u_s)$; for the other notations we refer to \eqref{eq5.2'}. Then, using BSDE \eqref{BSDE-P} and applying the It\^{o} formula to $p_s(P)(x_s^{u}-\hat{x}_s)$, we obtain, for $P\in\mathcal{P}_{\{\psi(\hat{u})\}}$,
\begin{equation}\label{equ.5.2b'}
\begin{array}{lll}
& &\displaystyle E_P\Big[ \hat{\Phi}_{x}(T)(x_{T}^{u}-\hat{x}_{T})+\int_0^T\hat{l}_{x}(t)(x_{t}^{u}-\hat{x}_{t})dt\Big]\\
&=&\displaystyle E_P\Big[ -\int_0^T(x_s^{u}-\hat{x}_s)\big(\hat{b}_x(s)p_s(P)ds+\hat{\sigma}_x(s)q_s(P)\langle B\rangle_s\big)\\
& &\quad+\displaystyle\int_0^T\big\{ p_s(P)\big(b^{u}(s)-\hat{b}(s)\big)ds+q_s(P)\big(\sigma^{u}(s)-\hat{\sigma}(s)\big)d\langle B\rangle_s\big\}\Big]\\
&=&\displaystyle E_P\Big[\int_0^T \big\{d(H^{u,P}(s)-\hat{H}^P(s))-(x^{u}_s-\hat{x}_s)d\hat{H}_x^P(s)\big\}\Big],
\end{array}
\end{equation}
and from the convexity of $H$ we conclude that
\begin{equation}\label{equ.5.2b}
\begin{array}{lll}
& &\displaystyle E_P\Big[ \hat{\Phi}_{x}(T)(x_{T}^{u}-\hat{x}_{T})+\int_0^T\hat{l}_{x}(t)(x_{t}^{u}-\hat{x}_{t})dt\Big]\\
&=&\displaystyle E_P\Big[\int_0^T \big\{d(H^{u,P}(s)-\hat{H}^P(s))-(x^{u}_s-\hat{x}_s)d\hat{H}_x^P(s)\big\}\Big]\\
&=&\displaystyle E_P\Big[\int_0^T \hskip -0.3cm\Big(\big\{d(H^{u,P}(s)\hskip -0.05cm -\hskip -0.05cm\hat{H}^P(s))\hskip -0.05cm -\hskip -0.05cm (x^{u}_s-\hat{x}_s)d\hat{H}_x^P(s)\hskip -0.05cm -\hskip -0.05cm (u_s-\hat{u}_s)d\hat{H}^P_v(s)\big\}\hskip -0.05cm +\hskip -0.05cm (u_s-\hat{u}_s)d\hat{H}^P_v(s)\Big)\Big]\\
&\ge &\displaystyle E_P\Big[\int_0^T(u_s-\hat{u}_s)d\hat{H}^P_v(s)\Big].
\end{array}
\end{equation}
Similarly we see that, for all $Q\in\mathcal{P}_{\{\varphi_4(\hat{x}_T)\}},$
\begin{equation}\label{equ.5.2ab}
\displaystyle E_Q\Big[\varphi'_4(\hat{x}_T)(x_{T}^{u}-\hat{x}_{T})\Big]\geq\displaystyle E_Q\Big[\int_0^T (u_s-\hat{u}_s)d\hat{H}_v^Q(s)\Big],
\end{equation}
for $\hat{H}^Q_v$ defined like $\hat{H}^P_v$, but with $(\tilde{p}_s(Q),\tilde{q}_s(Q))$ instead of $(p_s(P),q_s(P))$.
Similarly, for all $R=(R_t)\in\mathcal{R}_{\{\varphi_5(\hat{x})\}}$, we have, $dt$-a.e.,
\begin{equation}\label{equ.5.2c}
\displaystyle E_{R_t}\Big[\varphi'_5(\hat{x}_t)(x_{t}^{u}-\hat{x}_{t})\Big]\ge\displaystyle E_{R_t}\Big[\int_0^t (u_s-\hat{u}_s)d\hat{H}_v^{R_t}(s)\Big],
\end{equation}
where $\hat{H}_v^{R_t}(s),\, s\in[0,T]$, is defined by \eqref{equ.5.2Ham}, but with the solution $(p_s(t,R_t),q_s(t,R_t))$ of BSDE \eqref{BSDE-Rt} instead of that of BSDE \eqref{BSDE-P}. Consequently,
from \eqref{equ.5.2a}, \eqref{equ.5.2b}, \eqref{equ.5.2ab} and \eqref{equ.5.2c}, for all $(P,Q,R)\in\mathcal{P}\{\hat{u}\},$
\begin{equation}\label{equ.5.2d}
\begin{array}{lll}
& &J(u)-J(\hat{u})\\
&\ge &\displaystyle E_P\Big[ \hat{\Phi}_{x}(T)(x_{T}^{u}-\hat{x}_{T})+\int_0^T\big(\hat{l}_{x}(t)(x_{t}^{u}-\hat{x}_{t})+\hat{l}_v(t)(u_t-\hat{u}_t)\big)dt\Big]\\
& &\quad +\displaystyle E_P[\hat{\Phi}_y(T)]E_Q\Big[\varphi'_4(\hat{x}_T)(x_{T}^{u}-\hat{x}_{T})\Big]+\displaystyle \int_0^TE_P[\hat{l}_y(t)]E_{R_t}\Big[\varphi'_5(\hat{x}_t)(x_{t}^{u}-\hat{x}_{t})\Big]dt\\
&\ge&\displaystyle E_P\Big[\int_0^T \big\{(u_s-\hat{u}_s)d\hat{H}_v^P(s)+\hat{l}_v(s)(u_s-\hat{u}_s)ds\big\}\Big]\\
& &\quad +\displaystyle E_P[\hat{\Phi}_y(T)] E_Q\Big[\int_0^T (u_s-\hat{u}_s)d\hat{H}_v^Q(s)\Big]+\displaystyle \int_0^TE_P[\hat{l}_y(t)]\displaystyle E_{R_t}\Big[\int_0^t (u_s-\hat{u}_s)d\hat{H}_v^{R_t}(s)\Big]dt.
\end{array}
\end{equation}
Finally, recalling the notations introduced in Subsection 5.1.2, we see that the latter expression in \eqref{equ.5.2d} coincides with $\Theta[P,Q,R]$  (see \eqref{Theta1}), i.e., because of the arbitrariness of $(P,Q,R)\in\mathcal{P}\{\hat{u}\}$ in \eqref{equ.5.2d} we conclude that
$$J(u)-J(\hat{u})\ge \sup_{(P,Q,R)\in\mathcal{P}\{\hat{u}\}}\Theta[P,Q,R](u-\hat{u})\ge 0,$$
where the latter inequality comes from the assumption of our statement. This proves the optimality of the control $\hat{u}$.
\hfill$\square$

\noindent\textbf{Example 5.1}.  We consider the following  linear-quadratic control problem. The state equation is given by
$$
\left\{\begin{array}{l}
d x_t^u=(Ax_t^u+B u_t) d t+(Cx_t^u+D u_t) d B_t, \\
x(0)=x \in \mathbb{R},
\end{array}\right.
$$
where $u\in\mathcal{U}$ and $A, B, C, D$ are constants. We associate the cost functional
$$
J(u)=\frac{1}{2} \hat{\mathbb{E}}\left[\int_{0}^{T}( (x_{t}^{u})^2+u_t^2) d t+ (x_{T}^{u})^2+\hat{\mathbb{E}}[(x_{T}^{u})^2] \right],~u\in \mathcal{U}.
$$
The stochastic optimal control problem consists in  minimizing the cost functional over $\mathcal{U}$.

We remark that the running cost and the terminal cost in the cost functional $J(u)$ do not satisfy (A.2), but one checks rather easily that our arguments apply also here, as $x^{u}\in M_G^p(0, T),\ x_T^u\in L_G^p(\Omega)$, for all $p\geq 1$. We see in particular that the adjoint BSDEs \eqref{BSDE-P} and \eqref{BSDE-Q} take the form
\begin{equation}\label{equ.5.4a}
\begin{array}{lll}
\hskip -0.5cm dp_s(P)\hskip -0.3cm&=&\hskip -0.3cm (-Ap_s(P)_s+\hat{x}_s)ds-Cq_s(P)d\langle B\rangle_s+q_s(P)dB_s+dN_s(P),\, t\in[0,T],\, p_T(P)=\hat{x}_T,\\
\hskip -0.5cm d\tilde{p}_s(Q)\hskip -0.3cm&=&\hskip -0.3cm -A\tilde{p}_s(Q)_sds-C\tilde{q}_s(Q)d\langle B\rangle_s+\tilde{q}_s(Q)dB_s+d\tilde{N}_s(Q),\, t\in[0,T],\, \tilde{p}_T(Q)=\hat{x}_T,
\end{array}
\end{equation}
 respectively. We also remark that the solution $(p(t,R_t),q(t,R_t),N(t,R_t))$ is identically equal to zero, for all $t\in[0,T]$, since the running cost $l$ only depends on $(x^{u},u)$. So, with the notation $\displaystyle\psi(\hat{u}):=\frac{1}{2} \Big(\int_{0}^{T}\big( (\hat{x}_{t})^2+\hat{u}_t^2\big) d t+(\hat{x}_{T})^2+\hat{\mathbb{E}}[(\hat{x}_{T})^2]\Big),$
Theorem \ref{th5.6} says that there exists $(P^*,Q^*)\in\mathcal{P}_{\{\psi(\hat{u})\}}\times\mathcal{P}_{\{(\hat{x}_T)^2\}}$ such that, for all $u\in\mathcal{U}$,
\begin{equation}\label{5.5a}
\begin{array}{lll}
\hskip -0.25cm 0\hskip -0.25cm &\le &\hskip -0.25cm\displaystyle (u_s-\hat{u}_s)\hskip -0.05cm\Big(Bdsdp_s(P^*,Q^*)+
\hat{u}_sdsdP^*+Dd\langle B\rangle_sdq_s(P^*,Q^*)\Big)\\
\hskip -0.25cm&=&\hskip -0.25cm (u_s-\hat{u}_s)\hskip -0.05cm\Big(Bds(p_s(P^*)dP^*+\tilde{p}_s(Q)dQ^*)+
\hat{u}_sdsdP^*+Dd\langle B\rangle_s(q_s(P^*)dP^*+\tilde{q}_s(Q^*)dQ^*)\Big).
\end{array}
\end{equation}
On the other hand, we see that our example also satisfies the assumptions (A.3)-(A.4). Consequently, we have the following:
\begin{lemma} For our linear-quadratic control problem of Example 5.1 the condition \eqref{5.5a} is a necessary but also sufficient optimality condition for an admissible control $\hat{u}$.
\end{lemma}

\section{Appendix}
\subsection{Appendix 1. An extension of the result of Section 4}
Let us consider a  function $f:\mathcal{P}_2(\mathbb{R}^d)\rightarrow \mathbb{R}$ which is Lipschitz, i.e., there exist $C>0$ such that
\begin{equation}\label{6.1}
|f(\mu)-f(\mu')|\leq CW_2(\mu,\mu'),~\mu,\mu'\in\mathcal{P}_2(\mathbb{R}^d).
\end{equation}
We put $F_f(\xi):=\displaystyle{\sup_{P\in\mathcal{P}}}f(P_{\xi}),~\xi\in L_G^2(\Omega;\mathbb{R}^d).$\

\begin{remark}\rm From (\ref{6.1}) one sees immediately that $F_f:L_G^2(\Omega;\mathbb{R}^d)\rightarrow\mathbb{R}$ is Lipschitz.

Indeed, we have
$$
|F_f(\xi)-F_f(\eta)|\leq C\sup_{P\in\mathcal{P}}W_2(P_{\xi},P_{\eta})\leq C(\hat{\mathbb{E}}[|\xi-\eta|^2])^{\frac{1}{2}},~\xi,\eta\in L_G^2(\Omega;\mathbb{R}^d).
$$
\end{remark}
\begin{lemma}
Let $\xi\in L_G^2(\Omega;\mathbb{R}^d)$, with $\hat{\mathbb{E}}[|\xi|^2I_{\{|\xi|\geq N\}}]\rightarrow0~(N\rightarrow\infty)$. Then,
$$
\mathcal{P}_{\{\xi\}}^{f}:=\{P\in\mathcal{P}:f(P_{\xi})=F_f(\xi)\}\neq\emptyset.
$$
\end{lemma}
\noindent \emph{Proof.} As $F_f(\xi)=\displaystyle{\sup_{P\in\mathcal{P}}}f(P_{\xi})\leq f(\delta_0)+CW_2(\delta_0,P_{\xi})\leq f(\delta_0)+C(\hat{\mathbb{E}}[|\xi|^2])^{\frac{1}{2}}<+\infty$, where $\delta_0$ is the Dirac measure at $0\in \mathbb{R}^{d}$, there exists $(P^l)_{l\geq1}\subset\mathcal{P}$ such that $f(P_{\xi}^l)\uparrow F_f(\xi)$, as $l\rightarrow\infty$. But, since $\mathcal{P}$ is weakly compact, we can extract a subsequence $(P^{l'})_{l'\geq1}\subset(P^l)_{l\geq1}$, and find some $P\in\mathcal{P}$ such that $P^{l'}\rightharpoonup P$ (weak convergence), as $l'\rightarrow\infty$. Then
\begin{equation}\label{wc}
 E_{P^{l'}}[\theta]\rightarrow E_{P}[\theta], \text{~as~} l'\rightarrow\infty,\text{~for all~} \theta\in L_{ip}(\Omega).
\end{equation}
 Following the argument of the proof of Lemma \ref{lem4.1} in Section 4, it follows that
\begin{equation}\label{6.2}
 E_{P^{l'}}[\varphi(\xi)]\rightarrow E_{P}[\varphi(\xi)], \text{~for all~} \varphi\in C_b^1(\mathbb{R}^d).
\end{equation}
Indeed, $\varphi(\xi)\in L_G^2(\Omega)$, and for any $\delta>0$, there is $\theta\in L_{ip}(\Omega)$ such that
$$
\sup_{Q\in \mathcal{P}}E_Q[|\varphi(\xi)-\theta|^2]\leq\delta,
$$
and so (\ref{6.2}) follows from (\ref{wc}). But (\ref{6.2}) means that $P^{l'}_{\xi}\rightharpoonup P_{\xi}$, as $l'\rightarrow\infty$. As, on the other hand,
$$
\sup_{l'\geq1}E_{P^{l'}}[|\xi|^2I_{\{|\xi|\geq N\}}]\leq\hat{\mathbb{E}}[|\xi|^2I_{\{|\xi|\geq N\}}]\rightarrow0~(N\rightarrow\infty),
$$
Hence, $f(P_{\xi})=\displaystyle{\lim_{l'\rightarrow\infty}}f(P_{\xi}^{l'})=F_f(\xi)$, i.e., $P\in\mathcal{P}_{\{\xi\}}^{f}.$
\hfill$\square$

Similar to Section 4 we have
\begin{lemma}\label{lemma6.2}
Let $\xi, ~\eta\in L_G^2(\Omega;\mathbb{R}^d)$, with $\hat{\mathbb{E}}[(|\xi|^2+|\eta|^2)I_{\{|\xi|+|\eta|\geq N\}}]\rightarrow0~(N\rightarrow\infty)$. Then, for all $0<\varepsilon_l \downarrow0$, as $l\rightarrow\infty$, and any $P^{l}\in\mathcal{P}^f_{\{\xi+\varepsilon_l \eta\}}$, $l\geq1,$ we have:

{\rm i)} There exists a subsequence $(P^{l'})_{l' \geq 1}\subset(P^{l})_{l \geq 1}$
and $P \in \mathcal{P}$ such that $P^{l'} \rightharpoonup P$, as $l' \rightarrow \infty$.

{\rm ii)} If, for some $P \in \mathcal{P}$, $P^{l} \rightharpoonup P$, as $l \rightarrow \infty$, then $P \in \mathcal{P}_{\{\xi\}}^{f}$.

\end{lemma}
\noindent \emph{Proof.} For $0<\varepsilon_l \downarrow0$ $(l\uparrow\infty)$, let  $P^{l}\in\mathcal{P}^f_{\{\xi+\varepsilon_l \eta\}}$, $l\geq1.$ From the weak compactness of $\mathcal{P}$ it follows that there is a subsequence $(P^{l'})_{l'\geq1}\subset(P^l)_{l\geq1}$ and some  $P\in\mathcal{P}$ such that $P^{l'}\rightharpoonup P$, as $l'\rightarrow\infty$, i.e., for all $\theta\in L_{ip}(\Omega)$, $E_{P^{l'}}[\theta]\rightarrow E_{P}[\theta]$, as $l'\rightarrow\infty$, and as $(E_{P^{l'}}[\cdot])_{l'\geq1}$ is dominated by $\hat{\mathbb{E}}[\cdot]$, we have $E_{P^{l'}}[\zeta]\rightarrow E_{P}[\zeta]$, as $l'\rightarrow\infty$, for all $\zeta\in L_G^1(\Omega)$.

 Hence, for $\zeta=\varphi(\xi), ~\varphi \in C_{b}^{1}(\mathbb{R}^{d})$,
$$
\Big|E_{P^{l'}}[\varphi(\xi+\varepsilon_{l'}\eta)]-E_P[\varphi(\xi)]\Big|\leq\Big|E_{P^{l'}}[\varphi(\xi)]-E_P[\varphi(\xi)]\Big|+C_{\varphi}\varepsilon_{l'}\hat{\mathbb{E}}[|\eta|]\rightarrow0,~l'\rightarrow\infty.
$$
This combined with
$$\sup_{l'\geq1} E_{P^{l'}}[|\xi+\varepsilon_{l'}\eta|^{2}I_{\{|\xi+\varepsilon_{l'}\eta|\geq N\}}]\rightarrow0,~N\rightarrow\infty$$
(Recall the assumption on $\xi$ and on $\eta$) yields $W_{2}(P_{\xi+\varepsilon_{l'}\eta}^{l'}, P_{\xi}) \rightarrow 0 ~(l' \rightarrow \infty)$ . Consequently,
$$
F_f(\xi+\varepsilon_{l'}\eta)=f(P_{\xi+\varepsilon_{l'}\eta}^{l'})\rightarrow f(P_{\xi}), \text{~as~} l'\rightarrow\infty,
$$
since $P^{l'}\in \mathcal{P}_{\{\xi+\varepsilon_{l'} \eta\}}^f$, while, on the other hand,
$$\big|F_f(\xi+\varepsilon_{l'}\eta)-F_{f}(\xi)\big| \leq C(\hat{\mathbb{E}}[|\varepsilon_{l'}\eta|^{2}])^{\frac{1}{2}} \rightarrow 0, ~l' \rightarrow \infty.$$
It follows that $f(P_{\xi})=F_f(\xi)$, i.e., $P^{l'}\rightharpoonup P\in\mathcal{P}_{\{\xi\}}^{f}$.\hfill$\square$

From Lemma \ref{lemma6.2} we get
\begin{prop}
$\Gamma\big(\mathcal{P}_{\{\xi+\varepsilon \eta\}}^f, \mathcal{P}_{\{\xi\}}^{f}\big) \rightarrow 0$, as $0<\varepsilon\downarrow 0$.
\end{prop}
The proof is analogous to that of Proposition \ref{prop4.2}, and so we omit it here.

\medskip

Our objective is to study the (right- and left-) differentiability of $\lambda \rightarrow F_f(\xi+\lambda \eta)=\displaystyle{\sup_{P\in\mathcal{P}}} f(P_{\xi+\lambda\eta})$, for $\xi,~ \eta \in L^{2}(\Omega; \mathbb{R}^{d})$, with $\hat{\mathbb{E}}[(|\xi|^2+|\eta|^2)I_{\{|\xi|+|\eta|\geq N\}}]\rightarrow0$, as $N\rightarrow\infty$.

For this we suppose that $f\!\!: \mathcal{P}_{2}(\mathbb{R}^{d}) \rightarrow \mathbb{R}$ is differentiable in Lion's sense with Lipschitz continuous derivative $\partial_{\mu}f\!:\mathcal{P}_{2}(\mathbb{R}^{d}) \times\mathbb{R}^{d}\rightarrow \mathbb{R}^d$. Recall (see \cite{24}) that
$f: \mathcal{P}_{2}(\mathbb{R}^{d}) \rightarrow \mathbb{R}$ is differentiable, if there exists a continuous function $ \partial_{m}f: \mathcal{P}_{2}(\mathbb{R}^{d}) \times\mathbb{R}^{d}\rightarrow \mathbb{R}$ with $\partial_{m}f(\mu, \cdot): \mathbb{R}^{d} \rightarrow \mathbb{R}$ differentiable, for all $\mu \in \mathcal{P}_{2}(\mathbb{R}^{d})$, such that
$$\lim_{0<\lambda\downarrow0}\frac{f\big((1-\lambda) \mu+\lambda \mu'\big)-f(\mu)}{\lambda}=\int_{\mathbb{R}^d} \partial_{m}f(\mu, y)(\mu'-\mu)(d y), \text{~for all~} \mu, ~\mu'\in \mathcal{P}_{2}(\mathbb{R}^{d}),$$
and the derivative of $f$ w.r.t. the measure $\mu$ is defined by $\partial_{\mu}f(\mu, y):=\partial_{y}(\partial_{m} f)(\mu, y), \text{~for all~} \mu \in \mathcal{P}_{2}(\mathbb{R}^{d}),~ y\in\mathbb{R}^d.$

\noindent Note, for all $P\in\mathcal{P}$,
\begin{equation}\label{6.3}
f(P_{\xi+\varepsilon\eta})=f(P_{\xi})+\int_{0}^{1} \partial_{\lambda}[f(P_{\xi+\lambda\varepsilon\eta})]d\lambda
=f(P_{\xi})+\varepsilon E_{P}[(\partial_{\mu}f)(P_{\xi}, \xi)\eta]+\varepsilon R_{P}^{\varepsilon}, ~\varepsilon\geq0,
\end{equation}
where $$R_{P}^{\varepsilon}=\int_{0}^{1} E_{P}\big[\big((\partial_{\mu}f)(P_{\xi+\lambda\varepsilon\eta}, \xi+\lambda\varepsilon\eta)-(\partial_{\mu}f)(P_{\xi}, \xi)\big)\eta\big] d \lambda,$$
and
\begin{equation}\label{6.4}
|R_{P}^{\varepsilon}| \leq C \varepsilon \hat{\mathbb{E}}[|\eta|^{2}], ~\varepsilon \geq 0.
\end{equation}
Let us put
$$
\begin{aligned}
&G(\lambda):=\sup _{P \in \mathcal{P}}\big\{f(P_{\xi})+\lambda E_{P}[(\partial_{\mu} f)(P_{\xi}, \xi) \eta]\big\}, ~\lambda \in \mathbb{R} .
\end{aligned}
$$
Then, $G:\mathbb{R} \rightarrow \mathbb{R}$ is convex, and, so, in particular, there exists its right-derivative $G_{+}^{'}(0)$ at $\lambda=0$. On the other hand, from our above estimates it follows that
$$
|F_{f}(\xi+\varepsilon \eta)-G(\varepsilon)| \leq C \varepsilon^2 \hat{\mathbb{E}}[|\eta|^{2}], ~\varepsilon \geq 0.
$$
Hence,
$$
\bigg|\frac{F_f(\xi+\varepsilon\eta)-F_f(\xi)}{\varepsilon}-G_{+}^{'}(0)\bigg|
\leq\bigg|\frac{G(\varepsilon)-G(0)}{\varepsilon}-G_{+}^{'}(0)\bigg|+C \varepsilon \hat{\mathbb{E}}[|\eta|^{2}] \rightarrow 0, \text {~as~} 0<\varepsilon\downarrow 0,
$$
i.e., the right-derivative of  $\varepsilon \rightarrow F_{f}(\xi+\varepsilon \eta)$  at  $\varepsilon=0$ exists and
$$
\lim_{0<\varepsilon\downarrow 0}\frac{F_f(\xi+\varepsilon\eta)-F_f(\xi)}{\varepsilon}=G_{+}^{'}(0).
$$
\begin{prop}
Let $f: \mathcal{P}_{2}(\mathbb{R}^{d}) \rightarrow \mathbb{R}$ be differentiable,  with Lipschitz derivative $\partial_{\mu}f:\mathcal{P}_{2}(\mathbb{R}^{d}) \times\mathbb{R}^{d}\rightarrow \mathbb{R}^d$, and let
$\xi, ~\eta \in L^{2}_G(\Omega; \mathbb{R}^{d})$, with $\hat{\mathbb{E}}[(|\xi|^2+|\eta|^2)I_{\{|\xi|+|\eta|\geq N\}}]\rightarrow0~(N\rightarrow\infty)$. Then,\\

\rm i) $\displaystyle{\lim_{0<\varepsilon\downarrow 0}}\frac{F_f(\xi+\varepsilon\eta)-F_f(\xi)}{\varepsilon}=\displaystyle{\sup_{P\in\mathcal{P}^f_{\{\xi\}}}} E_{P}[(\partial_{\mu} f)(P_{\xi}, \xi) \eta];$

\rm ii)$\displaystyle{\lim_{0>\varepsilon\uparrow 0}}\frac{F_f(\xi+\varepsilon\eta)-F_f(\xi)}{\varepsilon}=-\displaystyle{\sup_{P\in\mathcal{P}^f_{\{\xi\}}}} \Big(-E_{P}[(\partial_{\mu} f)(P_{\xi}, \xi) \eta]\Big).$
\end{prop}
\noindent \emph{Proof.} We remark that ii) follows from i) by replacing in i) $\eta$ by $(-\eta)$. Let us prove i).

For this, using Lemma \ref{lemma6.2}, let $P^l \in \mathcal{P}^f_{\{\xi+\varepsilon_l \eta\}},~l \geq 1$, and $P\in\mathcal{P}$, such that, for $0 <\varepsilon_l\downarrow 0$ $(l \rightarrow \infty)$, $W_{2}(P_{\xi+\varepsilon_l \eta}^l, P_{\xi}) \rightarrow 0$. Then $P\in\mathcal{P}^f_{\{\xi\}}$. Thanks to (\ref{6.3}) and (\ref{6.4})
\begin{equation*}
\begin{aligned}
&\frac{F_f(\xi+\varepsilon_l\eta)-F_f(\xi)}{\varepsilon_l}\le\frac{f(P_{\xi+\varepsilon_l \eta}^l)-f(P_{\xi}^l)}{\varepsilon_l}=E_{P^l}[(\partial_{\mu} f)(P_{\xi}^l, \xi) \eta]+R_{P^l}^{\varepsilon_l}\\
 =&~E_{P^l}[(\partial_{\mu} f)(P_{\xi}^l, \xi) \eta]+O(\varepsilon_l).
\end{aligned}
\end{equation*}
Moreover, as $W_{2}(P_{\xi+\varepsilon_l \eta}^l, P^l_{\xi}) \leq \varepsilon_l \big(\hat{\mathbb{E}}[|\eta|^{2}]\big)^{\frac{1}{2}}\rightarrow 0~(l \rightarrow \infty)$, also $W_{2}(P_{\xi}^l, P_{\xi})\rightarrow 0$, and so $E_{P^l}[(\partial_{\mu} f)(P_{\xi}^l, \xi) \eta]\rightarrow E_{P}[(\partial_{\mu} f)(P_{\xi}, \xi) \eta],~l \rightarrow \infty.$ This shows that
$$
\varlimsup_{0<\varepsilon_{l}\downarrow 0} \frac{F_f(\xi+\varepsilon_l\eta)-F_f(\xi)}{\varepsilon_l} \leq \sup_{P\in \mathcal{P}_{\{\xi\}}^{f}}E_{P}[(\partial_{\mu} f)(P_{\xi}, \xi) \eta].
$$
On the other hand, for all $Q \in \mathcal{P}_{\{\xi\}}^{f}$,
$$
\frac{F_f(\xi+\varepsilon_l\eta)-F_f(\xi)}{\varepsilon_l}\geq
\frac{f(Q_{\xi+\varepsilon_l \eta})-f(Q_{\xi})}{\varepsilon_l}\rightarrow E_{Q}[(\partial_{\mu} f)(Q_{\xi}, \xi) \eta],~l \rightarrow \infty.
$$
This proves i),
\begin{equation}\label{6.5}
\lim_{0<\varepsilon\downarrow0}
\frac{F_f(\xi+\varepsilon\eta)-F_f(\xi)}{\varepsilon}=\sup_{Q\in\mathcal{P}^f_{\{\xi\}}} E_{Q}[(\partial_{\mu} f)(Q_{\xi}, \xi) \eta].
\end{equation}
\begin{flushright}
$\square$
\end{flushright}
\begin{remark}\rm Let $\xi,~ \eta \in L_{G}^{2}(\Omega ; \mathbb{R}^{d})$, such that $\hat{\mathbb{E}}[(|\xi|^2+|\eta|^2)I_{\{|\xi|+|\eta|\geq N\}}]\rightarrow0~(N\rightarrow\infty)$, $\varphi \in C_{b}^{1}(\mathbb{R}^{d}), ~ f(P_{\vartheta}):=E_{P}[\varphi(\vartheta)]$, and
$$
F_{f}(\vartheta):=\sup_{P\in\mathcal{P}} E_{P}[\varphi(\vartheta)], ~\vartheta\in L_{G}^{p}(\Omega ; \mathbb{R}^{d}).
$$
Then, as $\partial_{\mu} f(\mu, y)=\nabla \varphi(y), ~(\mu, y) \in \mathcal{P}_{2}(\mathbb{R}^{d}) \times \mathbb{R}^{d}$,
(\ref{6.5}) gives the result of Section 4 , but only for $\xi, ~\eta \in L_{G}^{2}(\Omega, \mathbb{R}^{d})$ with $\hat{\mathbb{E}}[(|\xi|^2+|\eta|^2)I_{\{|\xi|+|\eta|\geq N\}}]\rightarrow0~(N\rightarrow\infty)$, while in Section 4 we have considered $\xi, ~\eta \in L_{G}^{1}(\Omega, \mathbb{R}^{d})$.
\end{remark}

\subsection{Appendix 2. A measurable selection theorem}

Let $\xi=(\xi_t)$ and $\eta=(\eta_t)$ be in $M_G^2(0,T;\mathbb{R})$ such that the following assumptions are satisfied:

\smallskip

\noindent(\textbf{B.1}) $\hat{\mathbb{E}}\big[|\xi_t-\xi_s|^2+|\eta_t-\eta_s|^2\big]\le C|t-s|,\, \, t,\ s\in[0,T],$ for some constant $C\geq 0.$
\begin{remark}\label{remA.1}
\rm Recall from Lemma \ref{lem2.9} that, for $\phi,\psi:\mathbb{R}\rightarrow\mathbb{R}$ Lipschitz functions,  the processes $\xi=(\xi_t=\phi(\hat{x}_t))$ and $\eta=(\eta_t=\psi(\hat{x}_t))$ satisfy assumption (B.1), where $\hat{x}$ is the solution of SDE (5.1) with $u=\hat{u}$ optimal control.

We also observe that, for all $\xi=(\xi_t),\ \eta=(\eta_t)\in M_G^2(0,T;\mathbb{R})$
satisfying (B.1), the function $t\mapsto \hat{\mathbb{E}}_{\{\xi_t\}}[\eta_t]$ is Borel measurable. Indeed, from (B.1) it follows that, for all $\varepsilon>0$, the function $t\mapsto \hat{\mathbb{E}}[\xi_t+\varepsilon\eta_t]-\hat{\mathbb{E}}[\xi_t]$ is continuous and, hence, Borel measurable. Consequently,  Lemma \ref{lem4.4} shows that also
\begin{equation}\label{A.1}
\displaystyle t\mapsto\hat{\mathbb{E}}_{\{\xi_t\}}[\eta_t]=\lim_{0<\varepsilon\downarrow 0}\frac{1}{\varepsilon}\big(\hat{\mathbb{E}}[\xi_t+\varepsilon\eta_t]-\hat{\mathbb{E}}[\xi_t]\big),\quad t\in[0,T].
\end{equation}
is a Borel function.
\end{remark}
\begin{theorem}\label{thA.2}
 Assume that $\xi=(\xi_t),\ \eta=(\eta_t)\in M_G^2(0,T;\mathbb{R})$  satisfy (B.1).  Then the mapping
\begin{equation}\label{A.2}
[0,T]\ni t\mapsto \mathcal{P}_{\{\xi_t|\eta_t\}}:=\big\{R\in \mathcal{P}_{\{\xi_t\}}\ :\ E_R[\eta_t]=\hat{\mathbb{E}}_{\{\xi_t\}}[\eta_t]\big\}\subset\mathcal{P}
\end{equation}
is a weakly measurable set-valued function with non empty values which are compact subsets of $(\mathcal{P},d)$ (Recall that $d$ is the L\'evy-Prokhorov metric on $\mathcal{P}$).
\end{theorem}
\begin{remark}\label{remA.2}
\rm Recall that, if $(X,\mathcal{G})$ is a measurable space and $Y$ a topological space,  a set-valued function $G:X\ni x\mapsto G(x)\subset Y$ for which the values $G(x)$ are non empty, closed subsets of $Y$, is called \textit{weakly measurable} if, for all open subset $\mathcal{O}$ of $Y$, it  holds $\{x\in X\,:\, G(x)\cap\mathcal{O}\not=\emptyset\}\in\mathcal{G}$.
\end{remark}
\begin{theorem}\label{th'A.1}
 Assume that $\xi=(\xi_t),\ \eta=(\eta_t)\in M_G^2(0,T;\mathbb{R})$  satisfy (B.1).  Then the mapping
$[0,T]\ni t\mapsto \mathcal{P}_{\{\xi_t|\eta_t\}}\subset\mathcal{P}$ admits a $\mathcal{B}([0,T])-\mathcal{B}(\mathcal{P})$- measurable selection ($\mathcal{B}([0,T])$ and $\mathcal{B}(\mathcal{P})$ are the Borel $\sigma$-field over $[0,T]$ and $(\mathcal{P},d)$, respectively), i.e., there is a selection $R_t\in \mathcal{P}_{\{\xi_t|\eta_t\}},\, t\in[0,T],$ such that the mapping $t\mapsto R_t$ is $\mathcal{B}([0,T])-\mathcal{B}(\mathcal{P})$- measurable.
\end{theorem}
The proof of this theorem is an immediate consequence of Theorem \ref{thA.2} and the Kuratowski and Ryll-Nardzewski measurable selection theorem (cf. \cite{26}).  For the proof of Theorem \ref{thA.2} we need the following well-known auxiliary result:
\begin{lemma}\label{lemmaA.1}
For a given measurable space $(X,\mathcal{G})$  and a separable metric space $(Y,d)$  a set-valued function $G:X\ni x\mapsto G(x)\subset Y$ with non empty, closed values is weakly measurable if and only if, for every $y$ from a dense subset of $Y$, the function $X\ni x\mapsto d(G(x),y)$ is $\mathcal{G}-\mathcal{B}(Y)$-measurable ($d(G(x),y)$ is the distance of $y$ to $G(x)$ in $(Y,d)$ and $\mathcal{B}(Y)$ is the Borel-$\sigma$-field on $Y$).
\end{lemma}
We are now able to give the proof of Theorem \ref{thA.2}.
\begin{proof}
Let us begin with observing that $\mathcal{P}_{\{\xi_t|\eta_t\}}$ is non empty, for every $t\in[0,T]$. Indeed, $\mathcal{P}_{\{\xi_t\}}\not=\emptyset$ and $\hat{\mathbb{E}}_{\{\xi_t\}}[\eta_t]=\displaystyle\max_{R\in\mathcal{P}_{\{\xi_t\}}}E_R[\eta_t]$ (see Lemma \ref{lem4.4}). On the other hand, by writing $\mathcal{P}_{\{\xi_t|\eta_t\}}=\big\{R\in \mathcal{P}\ :\ E_R[\xi_t]= \hat{\mathbb{E}}[\xi_t],\,  E_R[\eta_t]= \hat{\mathbb{E}}_{\{\xi_t\}}[\eta_t]\big\},$ we see easily that $\mathcal{P}_{\{\xi_t|\eta_t\}}$ is closed and, hence, also compact, as $\mathcal{P}$ is.

Now, for any sequence $0<\varepsilon_k\searrow 0$ ($k\nearrow+\infty$), we consider

\centerline{$\displaystyle \mathcal{P}^k_{\{\xi_t|\eta_t\}}:=\big\{R\in \mathcal{P}\ :\ E_R[\xi_t]>  \hat{\mathbb{E}}[\xi_t] -\varepsilon_k,\,  E_R[\eta_t]> \hat{\mathbb{E}}_{\{\xi_t\}}[\eta_t]-\varepsilon_k\big\}\supset\mathcal{P}_{\{\xi_t|\eta_t\}},\, k\ge 1.$}

\noindent Because of the compactness of $(\mathcal{P},d)$ also the closure $\overline{\mathcal{P}^k_{\{\xi_t|\eta_t\}}}\subset(\mathcal{P},d)$ is compact, and

\centerline{$\displaystyle \mathcal{P}_{\{\xi_t|\eta_t\}}=\bigcap_{k\ge 1}\downarrow\overline{\mathcal{P}^k_{\{\xi_t|\eta_t\}}}.  $}

\noindent As $(\mathcal{P},d)$ is compact, this space is in particular separable, i.e., we can choose a dense countable subset $\mathcal{D}\subset\mathcal{P}.$ Let us put  $\mathcal{D}^k_{\{\xi_t|\eta_t\}}   =\mathcal{D}\cap\mathcal{P}^k_{\{\xi_t|\eta_t\}},\, t\in[0,T],\ k\ge 1.$ As for all $R\in \mathcal{P}^k_{\{\xi_t|\eta_t\}}$ there exists $(R_\ell)_{\ell\ge 1}\subset\mathcal{D}$ s.t. $R_\ell\rightharpoonup R$ and, thus, also $E_{R_\ell}[\xi_t]\rightarrow E_R[\xi_t]$ and $E_{R_\ell}[\eta_t]\rightarrow E_R[\eta_t]$ ($\ell\rightarrow +\infty$), it follows that $\mathcal{D}^k_{\{\xi_t|\eta_t\}}\subset\overline{\mathcal{P}^k_{\{\xi_t|\eta_t\}}}$ is dense.

For $k\ge 1$, put

\centerline{$\displaystyle F_k(t,Q):=\mbox{dist}_{\overline{\mathcal{P}^k_{\{\xi_t|\eta_t\}}}}(Q)\big(=\inf\big\{d(Q,R): R\in \overline{\mathcal{P}^k_{\{\xi_t|\eta_t\}}}\big\}\big),\, (t,Q)\in[0,T]\times\mathcal{P}.$}

\noindent Then,
\begin{equation*}\begin{array}{lll}
F_k(t,Q)&=&\displaystyle \inf\{d(R,Q): R\in \mathcal{D}^k_{\{\xi_t|\eta_t\}}\}\\
&=&\displaystyle \inf_{R\in\mathcal{D}}\chi(d(R,Q),\alpha_k(t,R),\beta_k(t,R)),
\end{array}\end{equation*}
where
\begin{equation*}\begin{array}{lll}
\alpha_k(t,R):&=&\big(E_R[\xi_t]-(\hat{\mathbb{E}}[\xi_t]-\varepsilon_k)\big)^-,\\
\beta_k(t,R):&=&\big(E_R[\eta_t]-(\hat{\mathbb{E}}[\eta_t]-\varepsilon_k)\big)^-,\\
\chi(\rho,\alpha,\beta):&=&\left\lbrace\begin{array}{lll}
&\rho,& (\alpha,\beta)=(0,0)\\
&+\infty,&  (\alpha,\beta)\not=(0,0).
\end{array}\right.
\end{array}\end{equation*}
Observe that, thanks to assumption (B.1), the functions $t\mapsto\alpha_k(t,R),\, \beta_k(t,R)$ are continuous and, hence, Borel measurable, for all $R\in\mathcal{D}$, and so $[0,T]\times\mathcal{P}\ni (t,Q)\mapsto \chi(d(Q,R),\alpha_k(t,R),$ $\beta_k(t,R))$ is Borel measurable (more precisely, $\mathcal{B}([0,T])\otimes\mathcal{B}(\mathcal{P}) - \mathcal{B}(\bar{\mathbb{R}})$-measurable),  for all $R\in\mathcal{D}$.  But as $\mathcal{D}$ is countable, also the infimum w.r.t. $R\in\mathcal{D}$ over these Borel functions is Borel measurable. Consequently, $F_k: [0,T]\times\mathcal{P}\ni (t,Q)\mapsto F_k(t,Q)$ is a Borel function, for all $k\ge 1.$

On the other hand, since $\displaystyle \mathcal{P}_{\{\xi_t|\eta_t\}}=\displaystyle\bigcap_{k\ge 1}\downarrow\overline{\mathcal{P}^k_{\{\xi_t|\eta_t\}}}$, we have
\begin{equation*}
F(t,Q):=\displaystyle\mbox{dist}_{\mathcal{P}_{\{\xi_t|\eta_t\}}}(Q)=\lim_{k\rightarrow+\infty}\uparrow\mbox{dist}_{\overline{\mathcal{P}^k_{\{\xi_t|\eta_t\}}}}(Q)=\lim_{k\rightarrow+\infty}\uparrow F_k(t,Q),\, (t,Q)\in[0,T]\times\mathcal{P},
\end{equation*}
and, hence, $F:[0,T]\times\mathcal{P}\rightarrow\mathbb{R}$ is Borel measurable.
From Lemma \ref{lemmaA.1} we get now the weak measurability of  the set-valued function $t\mapsto \mathcal{P}_{\{\xi_t|\eta_t\}}.$
\end{proof}


\end{document}